\documentclass{amsart}
\usepackage{amssymb,latexsym,graphics}
\usepackage[mathscr]{eucal}

\newtheorem{thm}{Theorem}[section]

\newtheorem{lem}[thm]{Lemma}
\newtheorem{prop}[thm]{Proposition}

\theoremstyle{remark}
\newtheorem{rem}[thm]{Remark}

\theoremstyle{definition}

\newcommand{\les}{{\lesssim}}

\newcommand{\lp}[2]{\Vert \, #1 \, \Vert_{#2}}
\newcommand{\snabla}{\slash \!\!\!\!\nabla}

\newcommand{\ret}{\vspace{.3cm}}

\begin{document}

\title[Angular Regularity and Strichartz Estimates]
{Angular Regularity and Strichartz Estimates for the Wave Equation}
\author{Jacob Sterbenz}
\address{Institute for Advanced Study and Princeton University, 
Princeton NJ, 08540}
\email{sterbenz@math.princeton.edu}
\author{appendix by Igor Rodnianski}
\address{Princeton University, 
Princeton NJ, 08540}
\email{irod@math.princeton.edu}
\thanks{The author is supported in part by an NSF postdoctoral fellowship.}
\subjclass{}
\keywords{}
\date{}
\dedicatory{}
\commby{}


\begin{abstract}
We prove here essentially sharp $L^q(L^r)$ linear and bilinear estimates 
for the  wave equations on Minkowski space where we assume the initial data
possesses additional regularity with respect to fractional powers of
the angular momentum operators $\Omega_{ij}:=
x^i\partial_j - x^j\partial_i$. In this
setting, the range of exponents  $(q,r)$ vastly improves over what is available
for the wave equations based on translation invariant derivatives of the 
initial data, or uniform decay of the solution.   
\end{abstract}

\maketitle

\section{Introduction: Classical Strichartz estimates, improvements 
for spherically symmetric data, and Knapp counterexamples}

The aim of this work is to prove mixed Lebesgue space estimates for
solutions to the linear wave equation on Minkowski space in a setting
where the initial data is assumed to possess extra regularity with
respect to weighted derivatives in the angular variable. These types of
estimates arise naturally in applications to the scale invariant  
global existence theory of non--linear wave equations which do not possess certain 
``null'' structures in their non--linearities. For example, in a companion to this
article, we use the estimates proved here to show global existence and scattering
for the (4+1) Yang--Mills equations in the Lorentz gauge for a certain
class of small, scale invariant initial data.\\

All of the estimates we prove
will be of ``Strichartz type'', i.e. $L^q(L^r)$ space--time
estimates for solutions of $\Box\phi = 0$. Due to the presence
of extra weighted angular regularity, we will get a significant gain over 
the usual estimates for the wave
equations which are based solely on translation 
invariant derivatives of the initial data.
What separates our estimates
from the ``classical'' Strichartz estimates, is that they 
are not solely based on the uniform decay of solutions to the wave equation.
Instead, we will exploit a certain ``wave packet'' 
structure these solutions exhibit in radial coordinates. This allows us to
decompose our waves into a sum of pulses, each of which remains coherent
for all time.\footnote{This fact makes the type of wave packets we use here
more simple to handle that those used in say \cite{W_bilin}. In fact, all of
our estimates are proved by directly integrating in time, without using any
complicated induction procedures or further analytical tools such as the 
$TT^*$ argument. However, there is a certain lack of coherence our
wave packets exhibit, which seems to limit their usefulness to situations where
one can use some angular regularity. This is because, strangely enough, our
wave packets become \emph{more} coherent for large times. For relatively
small times, close to the ``focusing'' time of our wave packets, much of
this coherence seems to get lost. See Remark \ref{coh_rem} after the
statement of Proposition \ref{psi_asym_prop}.}
These pulses can then be treated
on an essentially individual basis.  This works well to prove a certain 
weak\footnote{Not in the sense of real interpolation. By weak we
simply mean an estimate which contains logarithmic divergences.} 
endpoint estimate for the range we are considering. We then interpolate
our endpoint with the endpoint from \cite{KT_str} to obtain the full
sharp (up to an $\epsilon$) set of estimates.\\

We shall also consider multilinear type estimates which involve 
weighted angular 
regularity on one or any number of the factors. For this set 
of estimates, we again obtain a vast improvement over the ``classical''
multilinear  estimates for the wave equation (e.g. compared with \cite{KRT}). 
We shall state and prove these estimates in the context of T. Tao's dual scale
machine for generating multilinear estimates (see \cite{T_semi}). All of the
estimates we prove here are sharp, up to an arbitrarily small loss of 
angular regularity, when tested against Knapp counterexamples. We will explain
this in more detail in the sequel.\\
 
After the proof of the estimates we describe here was discovered by
the first author, a shorter proof was found by Igor Rodnianski for the case of
$3\leqslant n$ spatial dimensions. We have included this in an appendix to the
present work and have elected to retain the discussion of our original
proof in the main body of the paper because it includes the
development of machinery that is
interesting in its own right and is perhaps more
flexible. Furthermore, our proof gives a lot of detailed and
interesting information as to how wave propagation works on high
spherical harmonics. In particular, we provide what seems to be a sharp localization
of band limited Hankel transforms, the type of which has been studied
by previous authors (see \cite{FZ_Bessel}). It is likely that the general procedure we 
employ here which uses this type phase space localization is
applicable to other dispersive
phenomena (e.g. wave or Schr\"odinger equations) on spherically 
symmetric backgrounds. \\
 
We now begin with a brief discussion of the usual Strichartz 
estimates for the wave equation. The standard reference for this material
at this point is the paper \cite{KT_str}. 
Let  $u_1$ be  a unit frequency solution\footnote{For definitions of the various
objects we present here, including Fourier transforms and mixed
Lebesgue spaces see the next section.} to the wave equation on
Minkowski space. By this we mean a function $u_1$ such that:
\begin{align}
        \Box u_1 \ &= \ 0 \ ,   &\hbox{supp}\, \{\widehat{u_1}(0),
        \widehat{\partial_t u_1}(0)\} \ &\subseteq \
        \Big\{ |\xi| \ \Big| \ \frac{1}{2} < |\xi| < 2 \Big\} \ . \label{basic_wave}
\end{align}
Without loss of generality, we may assume that $u_1$ is of the form
$u_1(t) = e^{it\sqrt{-\Delta}}f_1$, for some unit frequency function 
$f_1$ of the spatial variable only.\footnote{For more detail,
see the notation list in Section \ref{not_sect}.} 
Then two of the most basic quantities
which determine the space--time behavior of $u_1$ are the following
estimates:
\begin{align}
        \lp{e^{it\sqrt{-\Delta}}f_1}{L^2_x} \ &= \ \ \ \ \ \ \ \ 
        \lp{f_1}{L^2_x} \ , 
        &\hbox{(Energy Estimate)} \label{basic_energy_est} \\
        \lp{e^{it\sqrt{-\Delta}}f_1}{L^\infty_x} \ &\lesssim \
        t^\frac{1-n}{2}  \ \lp{f_1}{L^1_x} \ . 
        &\hbox{(Dispersive Estimate)} \label{disp_est}
\end{align}
By interpolating between \eqref{basic_energy_est} and \eqref{disp_est}, 
and using some standard duality arguments (specifically the $TT^*$ method),
it is possible to show the following set of space--time estimates
for $u_1$:\\

\begin{thm}[``Classical'' Strichartz estimates including endpoints
(see \cite{KT_str}, \cite{GB}, and \cite{ST_classical})]
\label{str_th}
Let $3 \leqslant n$ be the number of spatial dimensions, and let
$\sigma=\frac{n-1}{2}$, then the following
estimate holds for $2\leqslant q$:
\begin{equation}
        \lp{e^{it\sqrt{-\Delta}}f_1}{L^q(L^r)} \ \lesssim \ 
        \lp{f_1}{L^2} \ , \label{str_est}   
\end{equation}
where $\frac{1}{q} + \frac{\sigma}{r} \leqslant \frac{\sigma}{2}$, with
the exception of the forbidden $L^2(L^\infty)$ endpoint on $\mathbb{R}^{3+1}$. 
\end{thm}\ret

\noindent A key facet of the proof of \eqref{str_est},
is that it does not rely
on any other property of the evolution operator $e^{it\sqrt{-\Delta}}$ than the
estimates \eqref{basic_energy_est} and \eqref{disp_est}. In other words,
the proof only uses the conservation of energy and the \emph{uniform}
decay estimate \eqref{disp_est}. Furthermore, the estimate \eqref{str_est} is sharp
in that one cannot improve the range of $L^q(L^r)$ indices stated above without
replacing the $L^2$ norm on the right hand side of \eqref{str_est}
with something else. This can be seen as follows: Let us consider initial data 
sets $\widehat{f^\epsilon_1} = \chi_{B^\epsilon}$, 
where $\chi_{B^\epsilon}$ is the indicator function of a radially directed
block of dimensions $1\times\epsilon\times\ldots\times\epsilon$, lying along
the $\xi_1$ axis between $1/2 < \xi_1 < 2$. Then a quick calculation using the
integral formula:
\begin{equation}
        e^{it\sqrt{-\Delta}} f^\epsilon_1\, (x) \ = \
        \int e^{2\pi i (t |\xi| + x\cdot\xi)} \chi_{B^\epsilon}(\xi)\, d\xi
        \ , \label{wave_osc_int}
\end{equation}
shows that one has $|e^{it\sqrt{-\Delta}} f^\epsilon_1\, (x)| 
\sim \epsilon^{n-1}$ on the space--time region $S^\epsilon_{t,x}$:
\begin{align}
        t \ &= O(\epsilon^{-2}) \ , \notag \\
        t + x_1 \ &= \ O(1) \ , \notag \\ 
        x' \ &= \ O(\epsilon^{-1}) \ . \notag
\end{align}
Based on this one sees that:
\begin{equation}
        \epsilon^{\frac{1}{2}(\frac{\sigma}{2} - \frac{1}{2} - 
        \frac{\sigma}{r})}\,
        \lp{f^\epsilon_1}{L^2} \ \sim  \ \epsilon^{n-1}\, 
        \lp{ \chi_{S^\epsilon_{t,x} }  }{L^2(L^r)} \ \lesssim \
        \lp{e^{it\sqrt{-\Delta}}f^\epsilon_1}{L^2(L^r)} \ , \label{knapp_eps}
\end{equation}
where $\sigma = \frac{n-1}{2}$. Therefore, the condition 
$\frac{1}{2} + \frac{\sigma}{r} \leqslant \frac{\sigma}{2}$ of Theorem
\ref{str_th} must be satisfied. We note here that initial data sets 
$f_1^\epsilon$ are commonly referred to as \emph{Knapp counterexamples}.\\

With the above consideration in mind, it is natural to wonder if that 
somehow the initial data $f_1$ were forced to be more evenly spread
out along the various radial directions in Fourier space, then one could 
gain some
improvement on the range of indices in \eqref{str_est}. However,
any such improvement must somehow involve another mechanism than just the 
estimates \eqref{basic_energy_est} and \eqref{disp_est}. This can be seen
simply from the fact that both of these estimates are sharp even for 
spherically symmetric data (e.g. looking at waves which focus at lager and 
larger times). Even so, it has been known
for some time that, with the help of additional arguments based
on the specific form of the integral representation for $e^{it\sqrt{-\Delta}}$, 
one can obtain a significant improvement over \eqref{str_est} for spherically
symmetric initial data (see for example \cite{SS_notes} and \cite{S_lnw}). 
To understand  how  this can happen, consider a unit frequency
radially symmetric initial data set $f_1^0$, such that $\widehat{f_1^0}$
is a smooth bump function of the radial variable $|\xi|$. Then by using
an integration by parts and stationary phase argument on the integral
representation \eqref{wave_osc_int}, in conjunction with the phenomena
of finite speed of propagation\footnote{The asymptotic 
\eqref{good_sph_asym} in the interior of the light--cone can be 
proved using integration by parts--stationary phase. In the exterior
of the light cone,
this asymptotic is true thanks to finite speed of propagation. 
Strictly speaking, the
functions $\{f^0_1(0), i\sqrt{-\Delta}f\, (0)\}$ are not compactly supported.
However, they are exponentially localized around the origin, so one may
recover \eqref{good_sph_asym} in the exterior via 
weighted energy estimates. We will give an independent proof of this
fact in a moment (see \eqref{psi_asym}) which does not rely on weighted energy
estimates.}
it is not difficult
to see that one has the asymptotics (because we can let derivatives fall on
$\widehat{f_1^0}$, of course this works for smooth enough Fourier data
even if it is not spherically symmetric):
\begin{equation}
        |e^{it\sqrt{-\Delta}}f^0_1\, (r)| \ \leqslant \ 
        \frac{C_M}{|t|^\frac{n-1}{2}( 1 + \big| |t| -r \big| )^M} 
        \ . \label{good_sph_asym}
\end{equation} 
That is, at time t, $e^{it\sqrt{-\Delta}}f^0_1$ 
is essentially the indicator function of an $O(1)$ spherical shell 
of radius $t$ multiplied 
by the amplitude $t^\frac{1-n}{2}$. Based on this, one can easily computes that:
\begin{equation}
        \lp{e^{it\sqrt{-\Delta}}f^0_1}{L^r_x} \ \lesssim \
        \frac{1}{(1 + |t|)^{(n-1)(\frac{1}{2} - \frac{1}{r})} } \ . \notag
\end{equation}
Therefore, in order for us to have that the $L^2(L^r)$ norm of 
$e^{it\sqrt{-\Delta}}f^0_1$ is finite, we need that 
$\frac{1}{2} < (n-1)(\frac{1}{2} - \frac{1}{r})$, or equivalently, that
$\frac{2(n-1)}{n-2} < r$. This is a vast improvement over the requirement
of $\frac{2(n-1)}{n-3} \leqslant r$ 
coming from the Knapp counterexamples. The key
to this improvement is that along with the uniform decay rate \eqref{disp_est},
the waves $e^{it\sqrt{-\Delta}}f^0_1$ are highly localized in physical space
along the radial variable.\\

Of course, an arbitrary spherically symmetric wave will not have the
localization \eqref{good_sph_asym}. However, it is possible to in a straightforward
manner chop up a unit frequency spherically symmetric wave into a sum of 
pieces, each of which satisfy a time translated version of the asymptotic
\eqref{good_sph_asym}. This is accomplished via a suitable physical space
localization
of the radially symmetric Fourier transform as follows:  We begin by rewriting
the integral formula \eqref{wave_osc_int} for radially symmetric
initial data $f_1$ as:
\begin{equation}
        e^{it\sqrt{-\Delta}}f_1\, (r) \ = \
        \frac{2\pi}{r^\frac{n-2}{2}}\int_0^\infty 
        e^{2\pi i t\rho}J_{\frac{n-1}{2}}(2\pi r\rho)
        \, \widehat{f_1}(\rho)\, \rho^\frac{n}{2} \ d\rho \ , 
        \label{sph_wave_osc_int}
\end{equation}
where $J_{\frac{n-1}{2}}(y)$ is the Bessel function of order $\frac{n-1}{2}$.
We now us the well known asymptotics for Bessel functions of relatively small
order (see \cite{Wbessel}):
\begin{equation}
      J_{\frac{n-1}{2}}(y) \ = 
      \begin{cases}
              \left( \frac{2}{\pi y} \right)^\frac{1}{2}\Big[
              \cos (y - \frac{n-2}{4}\pi)\cdot m_1(y) \ - \\
              \ \ \ \ \ \ \ \ \ \ \ \ \ \ \ \ \ \ 
              \sin (y - \frac{n-2}{4}\pi)\cdot m_2(y)\Big] \ , 
              &\hbox{ for $1 \leqslant y$ } \ , \\
              y^\frac{n-1}{2}\cdot m_3(y) \ , 
              &\hbox{ for $0 \leqslant y \leqslant 1$ } \ . 
      \end{cases}\label{bessel_asym}
\end{equation}
Here the function $m_3$ is $C^\infty$, and the remaining $m_i$ have
asymptotic expansions:
\begin{align}
        m_1(y) \ &= \ \sum_k C_{1,k}y^{-2k} \ , \notag \\
        m_2(y) \ &= \ \sum_k C_{2,k}y^{-2k-1} \ , \notag
\end{align}
as $y\to \infty$. In other words, the functions $m_1(2\pi r\rho)$ and 
$m_2(2\pi r\rho)$ are $C^\infty$ with derivatives
in $\rho$ uniformly bounded for all $\frac{1}{2} \leqslant \rho \leqslant 2$ and
$2 \leqslant r$. Substituting the asymptotic \eqref{bessel_asym} into the 
integral formula, we may assume without loss of generality that we are
trying to bound integrals of the form:
\begin{equation}
        I^\pm (t,r) \ = \ \frac{1}{(1 + r )^\frac{n-1}{2}}\int_{-\infty}^\infty
        e^{2\pi i (t \pm r)\rho}\, m^\pm (r,\rho)\, \chi_{(1/4,4)}(\rho)\, 
        \widehat{f_1}(\rho)\ d\rho \ , \label{basic_sph_int}
\end{equation}
where $m^\pm$ is a smooth function with derivatives in $\rho$ uniformly
bounded for all $0 \leqslant r $, and $\chi_{(1/4,4)}$ is a smooth bump
function on the interval $(\frac{1}{4},4)$. It is now apparent that the integrals
in \eqref{basic_sph_int} are essentially time translated inverse
Fourier transforms of a one dimensional unit frequency function. Therefore,
we can localize these integrals in physical space (i.e. the $t\pm r$ variable)
on an $O(1)$ scale. This can be accomplished with the help of the so called
$\varphi$--transform (see \cite{FJWlp}), which is just a smoothed out
redundant (over--sampled) version of the classical Shannon sampling. 
Since the function $\widehat{f_1}$ is compactly supported
in the interval $(0,4)$, we may take its Fourier series development:
\begin{align}
        \widehat{f_1}(\rho) \ = \ \sum_k \ c_k \, e^{i\frac{\pi}{2} k\rho} 
        \ \ \ \ ,&
        &\rho \ &\in \ (0,4) \ . \label{f_fourier_devel}
\end{align}
An important thing to notice here is that we can recover the $L^2$ of
$f_1$ as a function on $\mathbb{R}^n$ in terms of the $\{c_k\}$:
\begin{equation}
        \lp{f_1}{L^2_x}^2 \ \sim \ \sum_k\ |c_k|^2
        \ . \label{f_recover_est}
\end{equation}
This can be seen from the Plancherel theorem and the fact that 
$\widehat f_1$ is unit frequency, so the volume part of the integral
in Fourier space which comes from integrating over spheres is $O(1)$.
Sticking the series \eqref{f_fourier_devel} into the the integrals 
\eqref{basic_sph_int} yields:
\begin{equation}
        I^\pm (t,r) \ = \ \sum_k\
        \frac{c_k}{(1 + r )^\frac{n-1}{2}}\, \psi^\pm_k(t,r) \ ,
        \label{I_wavelet_exp}
\end{equation}
where:
\begin{equation}
        \psi^\pm_k(t,r) \ = \ \int_{-\infty}^\infty
        e^{2\pi i (t \pm r + \frac{k}{4})\rho}\, 
        m^\pm (r,\rho)\, \chi_{(1/4,4)}(\rho)\ d\rho \ . \notag
\end{equation}
Integrating by parts as many times as necessary in the above formula, we see
that we have the asymptotic:
\begin{equation}
        |\psi^\pm_k(t,r)| \ \leqslant \ 
        \frac{C_M}{(1 + | t \pm r + \frac{k}{4}|)^M} \ . \label{psi_asym}
\end{equation}
Using the expansion \eqref{I_wavelet_exp} and the asymptotic \eqref{psi_asym} 
we can directly compute that for $2\leqslant p$:
\begin{align}
        \lp{I^\pm (t,\cdot)}{L^p_x}^p \ &\lesssim \
        \int_{-\infty}^\infty \left( \sum_k\ c_k\, \psi^\pm_k(t,r) \right)^p
        (1 + r)^{(n-1)(1-\frac{p}{2})}\ dr \ , \notag \\
        &\lesssim \ \int_{-\infty}^\infty \left( \sum_k\ \frac{c_k}{
        (1 + | t \pm r + \frac{k}{4}|)^{3}}\right)^p
        (1 + r)^{(n-1)(1-\frac{p}{2})}\ dr \ , \notag \\
        &\lesssim \ \sum_k\
        \int_{-\infty}^\infty \ \frac{|c_k|^p}{
        (1 + \big| |t + \frac{k}{4} | - r \big|)^2}\,
        (1 + r)^{(n-1)(1-\frac{p}{2})}\ dr \ , \notag
\end{align}
The manipulation to get the last line above follows from H\"olders
inequality and the fact that $2\leqslant p$. By integrating each expression
in this line term by term and using the inclusion $\ell^2 \subseteq \ell^p$ for
$2\leqslant p$ we arrive at the bound: 
\begin{align}
        \lp{I^\pm (t,\cdot)}{L^p_x} \ &\lesssim \ \left(\sum_k\
        \frac{|c_k|^p}{( 1 + | t + \frac{k}{4} | )
        ^{p (n-1)(\frac{1}{2} - \frac{1}{p})} } \right)^\frac{1}{p} \ , 
        \notag \\
        &\lesssim \ \left(\sum_k\
        \frac{|c_k|^2}{( 1 + | t + \frac{k}{4} | )
        ^{2 (n-1)(\frac{1}{2} - \frac{1}{p})} } \right)^\frac{1}{2} \ . \notag
\end{align}
Testing this last expression for $L^2$ in time we see that:
\begin{equation}
       \lp{I^\pm}{L_t^2(L^p_x)}^2 \ \lesssim \
       \sum_k \ |c_k|^2 \, \int_{-\infty}^\infty \
       ( 1 + | t + \frac{k}{4} | )^{-2 (n-1)(\frac{1}{2} - \frac{1}{p})}\ dt
       \ . \notag \\
\end{equation}
Now as long as $1 < 2 (n-1)(\frac{1}{2} - \frac{1}{p})$, or equivalently
$\frac{2(n-1)}{n-2} < p$, we have the result:
\begin{equation}
        \lp{I^\pm}{L_t^2(L^p_x)}^2 \ \lesssim \ \sum_k \ |c_k|^2 \ . \notag
\end{equation}
Using the characterization \eqref{f_recover_est}, we have shown that:\\

\begin{prop}[Unit frequency Strichartz estimates for spherically symmetric data]
Let $u_1$ be a unit frequency, spherically symmetric solution to the equation
$\Box u_1 = 0$ in $3 \leqslant n$ (spatial) dimensions, then the following 
space--time estimates hold:
\begin{equation}
        \lp{u_1}{L^2(L^r)} \ \lesssim \ \lp{u_1(0)}{L^2} \ , 
        \label{str_sph_data}
\end{equation}
where $r$ satisfies the bound  $\frac{2(n-1)}{n-2} < r$.
\end{prop}\ret
\noindent Interpolating \eqref{str_sph_data} with the energy estimate
\eqref{basic_energy_est}, and by various rescalings and 
using Littlewood--Paley theory (as in the
work \cite{KT_str}, we will discuss this in more detail in the sequel), the 
estimate \eqref{str_sph_data} can be extended to general spherical initial data
and other $L^q(L^r)$ spaces in a straight forward manner. We record this as:\\

\begin{thm}[Strichartz estimates for spherically symmetric initial data]
\label{full_str_th_sph_data}
Let $u$ be a spherically symmetric function on $\mathbb{R}^{n+1}$ 
such that $\Box u = 0$, and set $\sigma_\Omega=n-1$, then 
the following estimates hold:
\begin{equation}
	\lp{u}{L^q(L^r)} \ \lesssim \ \left(\lp{u(0)}{\dot H^\gamma} +
	\lp{\partial_t u\, (0)}{\dot H^{\gamma-1}}\right) \ , 
	\label{full_sph_str_est}
\end{equation}
where $\frac{1}{q} + \frac{\sigma_\Omega}{r} < \frac{\sigma_\Omega}{2}$, and
$\frac{1}{q} + \frac{n}{r} = \frac{n}{2} - \gamma$.
\end{thm}\ret

\begin{rem}
In the language of \cite{KT_str}, Theorem \ref{full_str_th_sph_data} says
that the range of indices for Strichartz estimates for spherically symmetric
initial data are non-sharp $(n-1)$ admissible. From the point of view
of uniform decay and the general machinery of \cite{KT_str}, this is like 
saying that these solutions morally decay like $t^{1-n}$. Of course, solutions
to wave equation in general (if you look everywhere in Minkowski space
and take the supremum) only decay like $t^\frac{1-n}{2}$. The
point here is that if they are spherically symmetric,
they only do so on a relatively thin set.
\end{rem}\ret\ret

We would now like to prove a result like \eqref{str_sph_data} for non-spherical
initial data. By looking at the Knapp counterexample calculation 
\eqref{knapp_eps},
we see that in order for us to avoid a contradiction, we need to replace the 
term $\lp{f^\epsilon_1}{L^2}$ on the left hand side of that equation by 
something that is of the order $\epsilon^{-2(\frac{\sigma}{2} - \frac{1}{2} - 
\frac{\sigma}{r})}\,\lp{f^\epsilon_1}{L^2}$. Now, based on our experience with
spherical data, it is natural to conjecture that we only need to replace
the norm $\lp{f^\epsilon_1}{L^2}$ with something that incorporates angular
regularity.\footnote{The regularity in the radial directions does not effect the
optimal estimates one gets unless you ask for $L^q(L^r)$ estimates with 
$q < 2$. In order to get these type of estimates, it is necessary to 
incorporate some decay of the initial data, even if it is spherically symmetric.
This can be easily seen by taking a spherically symmetric wave $\psi$ with the 
asymptotic \eqref{good_sph_asym}, and then forming the wave $\Psi = \sum_{k=0}^N
\psi(t - k )$. Then at time zero one has $\lp{\Psi(0)}{L^2} \sim N^\frac{1}{2}$,
but one gets $\lp{\Psi}{L^q(L^r)} \sim N^\frac{1}{q}$ for 
$\frac{2(n-1)}{n-2} < r < \infty$. Therefore it is not possible to take 
$q < 2$.} What we will do then, is test for smoothness of $f^\epsilon_1$ in the
angular direction in Fourier space. This is done by using the 
infinitesimal generators of the rotations on Euclidean space:
\begin{equation}
        \Omega_{i,j} \ := \ x_i\partial_j - x_j\partial_i \ . 
        \label{rotation_fields}
\end{equation}
One sees immediately that:
\begin{equation}
        \lp{f^\epsilon_1}{|\Omega|^{-1}L^2}^2 \ := \
        \sum_{i<j} \ \lp{\Omega_{ij} f^\epsilon_1}{L^2}^2 \ = \
        \sum_{i<j} \ \lp{\Omega_{ij} \widehat{f^\epsilon_1}}{L^2}^2 \ \sim \
        \epsilon^{-2}\ \lp{f^\epsilon_1}{L^2}^2 \ . \label{contains_L2_omega}
\end{equation}  
Interpolating this with the identity, we see that:
\begin{equation}
        \lp{f^\epsilon_1}{|\Omega|^{-s} L^2} \ := \
        \lp{|\Omega|^s f^\epsilon_1}{L^2} \ \sim \ \epsilon^{-s} 
	\ , \label{Omega_of_f}
\end{equation}
where $|\Omega|^s = (-\Delta_{sph})^\frac{s}{2}$, and:
\begin{equation}
        \Delta_{sph} \ := \ \sum_{i<j} \Omega_{ij}^2 \ , \label{sph_lap}
\end{equation}
is the Laplacian on the sphere
of radius $r$.\footnote{Fractional powers of the operator $-\Delta_{sph}$
can be defined in the usual way via spectral resolution. Furthermore,
the interpolation identity: $(L^2,|\Omega|^{-1} L^2)_t = |\Omega|^{-t} L^2$
can easily be shown using spectral resolution and interpolation of weighted
$\ell^2$ sequence spaces. We'll discuss this in more detail in just a bit.}
Therefore, it is natural to expect that if we add 
$s=2(\frac{\sigma}{2} - \frac{1}{2} - \frac{\sigma}{r})$ angular
derivatives to the left hand side of \eqref{knapp_eps}, we may get a true
estimate in the range $\frac{2(n-1)}{n-2} < r$. Furthermore, an estimate of
this kind would be sharp. This is precisely what we will prove in
dimensions $4\leqslant n$, with an 
$\epsilon$ loss of angular regularity:\\

\begin{thm}[Strichartz estimates for angularly regular 
data]\label{full_ang_str_th}
Let $4 \leqslant n$ be the number of spatial dimensions, and let $u$ be a
solution to the homogeneous wave equation $\Box u = 0$.  Let 
$\sigma_\Omega=n-1$ be the angular wave admissible Strichartz
exponent, and let $\sigma=\frac{n-1}{2}$ be the classical 
wave admissible Strichartz exponent. Then 
for every $0 < \epsilon$, there is a $C_\epsilon$ depending only on $\epsilon$
such that the following set of estimates hold:
\begin{equation}
        \lp{u}{L^q(L^r)} \ \lesssim \ C_\epsilon\,\left(
        \lp{  \langle\Omega\rangle^s
        u(0)}{\dot H^\gamma} + \lp{  \langle\Omega\rangle^s
        \partial_t u\, (0)}{\dot H^{\gamma-1}}\right)
	\ , \label{full_ang_str_est}   
\end{equation}
where we have that $r\neq\infty$, \ 
$s= (1 + \epsilon)(\frac{n-1}{r} + \frac{2}{q}- \frac{n-1}{2})$,
\ $\frac{1}{q} + \frac{n}{r} =  \frac{n}{2} -\gamma$\ ,
\ $\frac{1}{q} + \frac{\sigma}{r} \geqslant \frac{\sigma}{2}$\ , and \
\ $\frac{1}{q} + \frac{\sigma_\Omega}{r} <  \frac{\sigma_\Omega}{2}$\ .
All of the
implicit constants in the above inequality depend on $n$, $q$, and $r$. Here:
\begin{equation}
	\lp{\langle\Omega\rangle^s u}{\dot H^\gamma}
	\ = \ \lp{u(0)}{\dot H^\gamma} + \lp{|\Omega|^s 
	u(0)}{\dot H^\gamma}\ , \notag
\end{equation}
with the analogous norm defined for $\lp{  \langle\Omega\rangle^s
\partial_t u\, (0)}{\dot H^{\gamma-1}}$.
\end{thm}\ret

\begin{rem}
A short calculation like the one done above shows that in fact
(modulo $\epsilon$ angular derivatives), all of the estimates 
\eqref{full_ang_str_est} are sharp when tested on Knapp counterexamples.
Therefore, in this sense, they are all endpoint estimates. It would
be interesting to try and remove the extra $\epsilon$ in these. 
We will not pursue this issue here, although we will do some extra
work in the sequel to recover a sharp $L^2$ dispersive estimate
which could be a start in this direction (see Proposition 
\ref{L2_disp_est_prop}). 
\end{rem}\ret

\begin{rem}
For the case of $n=2,3$ spatial dimensions, we will also prove an
estimate of the type $L^{2+}(L^\infty)$ and $L^2(L^{4+})$ respectively
which involves $\frac{1}{2}$ an
angular derivative. However, to obtain the full range of
\eqref{full_ang_str_est} in the case of $n=3$ spatial dimensions would require one
to prove an $L^2(L^\infty)$ Strichartz estimate that involves
$\epsilon$ angular derivatives. While it seems that this type of estimate
is out of the reach of methods we use here, it should be attainable
using the recent method wave packet of Wolff \cite{W_bilin}. In fact, all of
the estimates  \eqref{full_ang_str_est} should be able to be proved
directly using that method, based on the fact that they correspond
to estimates that should be true for angularly separated initial data
with no extra angular regularity.
\end{rem}\ret\ret

The remainder of this paper is laid out as follows. In the next section, we list 
briefly some of the basic notations we use here. \\

In the third section, we list some standard facts about analysis on
the sphere that will be useful in the sequel, including formulas for
Hankel transforms and Littlewood--Paley--Stein theory on the sphere. We then
use this machinery to reduce the proof of Theorem \ref{full_ang_str_th}
to a suitable set of ``endpoint'' estimates.\\

In the fourth section, we discuss the main tool to be used in this
paper: a $\varphi$--type transform for the Hankel transform. This leads
us to consider the localization properties in physical space of the bandwidth
limited Hankel transform. In particular, we provide detailed
asymptotics for our  Hankel--$\varphi$ transform which
will form the backbone the Strichartz estimates to be proved here.\\
   
In the fifth, we prove an $L^2$ dispersive estimate for
the wave equation based on angular regularity with respect to the 
momentum operators \eqref{rotation_fields}. This can be 
interpolated with the energy estimate to prove our ``endpoints'' directly
via integration in time, 
avoiding to use of any other analytic machinery such as usual $TT^*$ process
or induction on scales. \\

In the sixth, we give a brief description of
how our dispersive estimate can be modified in a straight forward
manner to accommodate multilinear phenomena. This shows one strength of the
method used here in that one gets the expected range of improved
multilinear estimates virtually for free out of the machinery developed.
We will also discuss why these multilinear estimates are sharp, by testing them on a 
multilinear analog of the Knapp counterexamples introduced above.\\

In an appendix to this paper, we provide a proof of the angularly
regular endpoints in the $3\leqslant n$ regime based on an idea of
Igor Rodnianski. The proof there is essentially independent of the
machinery we develop here (modulo a somewhat similar setup in terms
of multilinear estimates) and instead relies on a direct calculation
involving the energy-momentum tensor for the wave equation.

\ret

\section{Basic notation}\label{not_sect}

We list here some of the basic notation used throughout this paper.
For quantities $A$ and $B$, we denote by $A \lesssim B$ to mean
that $A \leqslant C\cdot B$ for some large constant $C$. The constant
$C$ may change from line to line, but will always remain
fixed for any given instance where this notation appears. We will also
use the notation $A\sim B$ if there exists a constant $C$ such that
$\frac{1}{C} \cdot A \leqslant B \leqslant C\cdot A$.\\

For a given function of tow variables, say $u(t,x)$, we denote the mixed
Lebesgue spaces norms $L^q(L^r)$ of $u$ via the formulas:
\begin{equation}
	\lp{u}{L^q(L^r)}^q \ := \ \int \ \lp{u(t)}{L_x^r}^q\ dt \ . \notag  
\end{equation}

For a given function of the spatial variable only, we denote its
Fourier transform as:
\begin{equation}
	\widehat{f}(\xi) \ := \ \int e^{-2\pi i \xi\cdot x} \, f(x) \ dx
        \ . \notag \\
\end{equation}
With this normalization of the Fourier variables, the Plancherel theorem
becomes $\lp{f}{L^2}=\lp{\widehat f}{L^2}$, and one has the Fourier
inversion formula:
\begin{equation}
	f(x) \ = \ \int e^{2\pi i \xi\cdot x} \, \widehat f(\xi) \ d\xi
        \ . \notag \\
\end{equation}
Using the Fourier transform, we define the homogeneous Sobolev
space of order $\gamma$ via the identity $\lp{f}{\dot H^\gamma} :=
\lp{|\xi|^\gamma \widehat f}{L^2}$.\\

For every integer $k$ we define the spatial Littlewood--Paley cutoff operator
by the formula:
\begin{equation}
	\widehat{P_k f} \ := \ p_k \, \widehat{f} \ , \label{lp_cut} 
\end{equation}
where $p_k(\xi)=p_0(2^{-k}\xi)$ and $p_0$ is a positive smooth bump function,
$p_0\equiv 1$ on the interval $(1,2)$ and zero off the interval
$(\frac{1}{2},4)$. With this notation we have the following consequence of the
Littlewood--Paley theorem for $2\leqslant q,r$ and $r < \infty$:
\begin{equation}
	\lp{u}{L^q(L^r)}^2 \ \lesssim \
	\sum_k \ \lp{P_k\, u}{L^q(L^r)}^2 \ . \label{sp_lp_sp_sum}
\end{equation}\\

For a given function of the spatial variable only, we denote its
forward and backward wave propagation via the formulas:
\begin{equation}
	W^\pm f\, (t,x) \ = \ e^{\pm i t \sqrt{-\Delta}}f\, (x) \ = \ 
	\int \ e^{\pm 2\pi i t |\xi|}e^{2\pi i x\cdot\xi}
	\, \widehat f(\xi)\ d\xi \ . \notag
\end{equation}
If $u$ is an arbitrary solution to the homogeneous wave equation,
$\Box u = 0$, with initial data $u(0) = f$ and $\partial_t u(0) = g$,
we may decompose it into forward and backward wave propagation in the 
following way:
\begin{equation}
	u(t,x) \ = \ \frac{1}{2i \sqrt{-\Delta}}\left(W^+g - W^-g\right) +
	\frac{1}{2}\left(W^+f + W^-f\right) \ . \label{wave_prop_decomp}
\end{equation}
Notice that in the case where $u=u_1$ is unit frequency, the functions:
\begin{align}
        h_1^+ \ &= \ \frac{1}{2i \sqrt{-\Delta}}g_1 + \frac{1}{2}f_1\ ,& 
	h_1^- \ &= \ -\frac{1}{2i \sqrt{-\Delta}}g_1 + \frac{1}{2}f_1\ , \notag
\end{align}
have $L^2$ norm comparable to the $\dot H^\gamma\times \dot H^{\gamma-1}$
norm of $(f_1,g_1)$. Therefore, in the sequel, we will always assume
that our unit frequency waves are of the form $e^{\pm it\sqrt{-\Delta}}h_1^\pm$,
and we replace:
\begin{equation}
        \lp{u_1(0)}{L^2} \ \sim \ (\lp{u_1(0)}{\dot H^\gamma} + 
	\lp{\partial_t u_1\, (0)}{\dot H^{\gamma-1}}) \ . \notag
\end{equation}

\ret

\section{Some results from analysis on the sphere.}

We list here some basic results from Fourier analysis is spherical coordinates
which will be used in our proof of Theorem \ref{unit_ang_str_th}. We have
already introduced the two basic differential elements of analysis on the
sphere, the rotation vector fields $\{\Omega_{ij}\}$ and the 
spherical Laplacian $\Delta_{sph}$. For every integer $0 \leqslant l$, there 
exists a finite dimensional set of functions $\mathcal{Y}_l$ on the sphere 
$\mathbb{S}^{n-1}\subseteq \mathbb{R}^n$ with satisfy the equation:
\begin{align}
        -\Delta_{sph} Y^l \ = \ l(n + l - 2)\, Y^l \ ,& &Y^l\in \mathcal{Y}_l
        \ . \notag
\end{align}
These sets of functions exhaust the set of eigenfunctions of $-\Delta_{sph}$
and can in fact be identified with the homogeneous polynomials $P^l$
on $\mathbb{R}^n$ of degree $l$ which satisfy:
\begin{equation}
        \big(\partial_r^2 + \frac{n-1}{r}\partial_r + \frac{1}{r^2} 
        \Delta_{sph}\big) P^l \ = \ 0 \ . \notag
\end{equation}
Setting $\omega_{n-1} := |\mathbb{S}^{n-1}|$, and introducing the natural
inner product on $\mathbb{S}^{n-2}$:
\begin{equation}
        \langle F , G \rangle \ := \ \omega_{n-1}^{-1} \int_{\mathbb{S}^{n-1}} 
        F(\omega)\, \overline{H}(\omega) \ d\omega \ , \notag
\end{equation}
we have the following basic properties of the $Y^l$:\\

\begin{lem}[Basic properties of spherical harmonics (see e.g. \cite{SWeuclid}]
\label{basic_sph_lem}
\ \ \ 
\begin{enumerate}
        \item  The dimension of the space $\mathcal{Y}_l $ is
        $|\mathcal{Y}^l| = \frac{1}{l}(n+ 2l - 2)\binom{n+l-3}{l-1}$. 
        \item  The spaces $\mathcal{Y}_l$ are mutually orthogonal. That
        is $\langle Y^l , Y^k \rangle = 0$ for $l\neq k$.
        \item  Let $\{Y^l_i\}_{i=1}^{|\mathcal{Y}_l|}$ be any orthonormal
        basis of $\mathcal{Y}_l$, then one has the following identity
        for all $\omega \in \mathbb{S}^{n-1}$: \ \ 
        $\sum_i \ |Y^l_i(\omega)|^2 \ = \ |\mathcal{Y}_l|$.
        \item  For each $Y^l\in \mathcal{Y}_l$, we have the identity:\\
        \ \ $\sum_{i < j} \lp{\Omega_{ij} Y^l}{L^2(\mathbb{S}^{n-1})}^2
        = l(n + l - 2)\lp{Y^l}{L^2(\mathbb{S}^{n-1})}^2 $
\end{enumerate}
\end{lem}\ret

>From now on, we fix an orthonormal basis $\{Y^l_i\}$ for each $\mathcal{Y}_l$.
For a given function $F\in L^2(\mathbb{S}^{n-1})$, we may expand it in the
$L^2$ sense along this basis as follows:
\begin{equation}
        F \ = \ \sum_{l,i} \ c_i^l\, Y^l_i \ . \label{sph_expansion}
\end{equation}
Using \eqref{sph_expansion}, we can define the action of $|\Omega|^s$
on this $F$ as follows:
\begin{equation}
        |\Omega|^s F \ = \sum_{l,i}
        \ [l(n+l-2)]^{\frac{s}{2}}\, c^l_i\, Y_i^l \ . \label{fractional_power}
\end{equation}
Then using item (2) and (4) of Lemma \eqref{basic_sph_lem}, and the 
fact that $\Omega_{ij} Y^l \in \mathcal{Y}_l$ ($\Omega_{ij} Y^l$
is a homogeneous harmonic polynomial of degree $l$ on $\mathbb{R}^n$),
we see that we have the identity:
\begin{equation}
        \sum_{i < j} \ \lp{\Omega_{ij} F}{L^2(\mathbb{S}^{n-1})}^2 \ = \
        \sum_{l,i} \ l(n+l-2) \,|c_i^l|^2 \ = \ \lp{|\Omega|\, F}
        {L^2(\mathbb{S}^{n-1})}^2 \ . \notag
\end{equation}
Using this, we see that there is, for every test function $f$ on 
$\mathbb{R}^n$ an equivalence of norms:
\begin{equation}
        \lp{f}{|\Omega|^{-1} L^2} \ = \ \lp{|\Omega| \, f}{L^2} \ , 
        \label{L2_Omega_norm}
\end{equation}
where $|\Omega|^{-1} L^2$ is the norm from line \eqref{contains_L2_omega}.
We also make the definition:
\begin{align}
	\lp{f}{H^s_\Omega}^2 \ &:= \ \lp{\langle \Omega \rangle^s f}{L^2}^2
	\ , \notag\\
	&:= \ \lp{f}{L^2}^2 + \lp{|\Omega|^s f}{L^2}^2 \ . \notag
\end{align}
Now, using the fact that $H^s_\Omega$ is of the form
$L^2(\ell^2_s)$, we have the following standard interpolation result 
(see \cite{BLinterp}):\\

\begin{prop}[Interpolation of the angular Sobolev
spaces $H^s_\Omega$]
\label{L2_Omega_interp_prop}
For any set or real numbers $s_1$ and $s_2$, we have the following 
interpolation spaces for $0 < t < 1$:
\begin{equation}
        \left( \ H^{s_1}_\Omega \ , \ H^{s_2}_\Omega\ \right)_t \ = \ 
        H^s_\Omega\ , \label{L2_Omega_interp} 
\end{equation}
where $s = (1-t)s_1 + t s_2$.
\end{prop}\ret

We will also use here some Littlewood--Paley theory in the 
angular variable\footnote{Notice that strictly speaking it will not
be necessary for us to use the Littlewood--Paley \emph{theorem} in the
angular variable due to our allowed loss of angular regularity. However, the use of
Littlewood--Paley \emph{cutoffs} in the angular variable will be essential for us.}.
We proceed in analogy with \eqref{lp_cut} and let $\theta_0$ be a smooth
bump function such that $\theta_0\equiv 1$ on the interval $(1,2)$ and
vanishing off the interval $(\frac{1}{2},4)$. For each $j\in\mathbb{Z}$
we denote its dyadic rescaling as $\theta_j(l) := \theta_0(2^{-j}l)$.
Defining $N=2^j$ and using the decomposition formula 
\eqref{sph_expansion}, we
define the angular frequency dyadic projections of a function $F$ on the 
sphere as:
\begin{equation}
	F_N \ := \ \sum_{l,i} \ c_i^l\, \theta_j(l)\, Y^l_i \ . \label{sph_lp_cut_def}
\end{equation}
We define $F_0$ to be the constant $c^0$, which is the average of $F$
over the sphere.
A formula similar to \eqref{sph_lp_cut_def}
can be used to define $f_N$ for a function on the whole
of $\mathbb{R}^n$.\\

For a given function $F$ on the sphere, we use line (1) from
\ref{basic_sph_lem}, as well as the fact that 
$\frac{1}{l}(n+ 2l - 2)\binom{n+l-3}{l-1} \sim N^\frac{n-2}{2}$
when $l\sim N$ to prove \emph{Bernstein's inequality} for the sphere:
\begin{align}
        | F_N | \ &\lesssim \ \sum_{l,i \ :\ l\sim N} \ |c_i^l|\cdot |Y^l_i|
	\ , \notag \\
	&\lesssim \ \sum_{l\sim N} (\sum_{i} |c_i^l|^2)^\frac{1}{2}\cdot
	N^\frac{n-2}{2} \ , \notag\\
	&\lesssim \ (\sum_{l,i \ :\ l\sim N} |c_i^l|^2)^\frac{1}{2}\cdot
	N^\frac{n-1}{2} \ , \notag\\
	&= \ N^\frac{n-1}{2}\, \lp{F}{L^2(\mathbb{S}^{n-1})}
	\ . \label{sph_bernstein}
\end{align}
By rescaling \eqref{sph_bernstein} to spheres of various radii, we have the following
result on all of $\mathbb{R}^n$ (for $N=0$, we of course replace the 
$N^\frac{n-1}{2}$ on the right hand side by $1$):
\begin{equation}
        \lp{f_N}{L^2_r(L^\infty(\mathbb{S}_r^{n-1}))} \ \lesssim \ 
	\, N^\frac{n-1}{2}\, \lp{r^{-\frac{n-1}{2}}f}{L^2(\mathbb{R}^n)} 
	\ . \label{Rn_sph_bernstein}
\end{equation}
Here $\mathbb{S}_r^{n-1}$ denotes the sphere of radius $r$ centered at the origin.\\

Next, we record the basic result which allows us to
generate certain square function expressions in terms of
the $f_N$ for $L^r$ spaces when $2\leqslant r$:\\

\begin{prop}[Littlewood--Paley--Stein theorem for the sphere (see \cite{Ssi}, 
\cite{Slp}, and \cite{STRlp})]\label{sph_lp_prop}
Let $\{\theta_j\}_{j=0}^\infty$ be any set of smooth functions such that there 
exists a $0 < \delta$ and $ \delta < \theta_j(l) < \frac{1}{\delta}$ for 
$l\in (2^j , 2^{j+1})$, and $\theta_j(l) = 0$ for $l\notin (2^{j-1}, 2^{j+1})$. 
Furthermore, let $\theta_j$ satisfy the bounds:
\begin{equation}
        \big|\frac{d^M}{dl^M} \theta_j\big| \ \leqslant \ C_M \cdot l^{-M} 
        \ . \notag
\end{equation}
Then one has that for any test function $F$ on $\mathbb{S}^{n-1}$, the following
bound on the ratio of norms holds for $1 < p < \infty$:
\begin{equation}
        \frac{1}{C_{p,\theta}} \ \ < \ \ 
        \lp{F}{L^p(\mathbb{S}^{n-1})} \  \Big/ \ \lp{ \Big(\sum_{j=0}^\infty \ 
        \big|
        \sum_{l,i} \ \theta_j(l)\, c^l_i\, Y^l_i \big|^2 \Big)^\frac{1}{2} }
        {L^p(\mathbb{S}^{n-1})}  \ \ <  \ \ 
        C_{p,\theta} \ . \label{sph_littlewood_paley}
\end{equation}
Here the constant $C_{p,\delta}$ depends only on $p$ and the $\{\theta_j\}$.
\end{prop}\ret

\noindent Using Proposition \ref{sph_lp_prop} along with the decomposition
\eqref{sp_lp_sp_sum}, we have the following estimates for
functions on space--time, for $2\leqslant q,r$ and $r< \infty$:
\begin{equation}
	\lp{u}{L^q(L^r)}^2 \ \lesssim \ \sum_{\substack{k\in\mathbb{Z}\\
	N\in 2^\mathbb{N}\cup \{0\}}}\ \lp{P_k\, u_N}{L^q(L^r)}^2 \ .
	\label{ang_lp_sp_sum}
\end{equation}\ret


\subsection{Reduction of Theorem \ref{unit_ang_str_th}  to a frequency localized
endpoint}
We now use the setup we have introduced above to reduce the proof of
Theorem \ref{unit_ang_str_th} to the following frequency localized
``endpoint'' estimate:\\

\begin{prop}[Endpoint unit frequency Strichartz estimate for angularly regular
data]\label{unit_ang_str_th}
Let $3 \leqslant n$ be the number of spatial dimensions, and let $u_{1,N}$ be a
unit frequency, angular frequency localized
solution to the homogeneous wave equation $\Box u_{1,N} = 0$.  Then
for every $0 < \eta$, there exists an $\frac{2(n-1)}{n-2} < r_\eta$, 
such that $r_\eta\to\frac{2(n-1)}{n-2}$ as $\eta\to 0$ and 
such that the following estimate holds:
\begin{equation}
        \lp{u_{1,N}}{L^2(L^{r_\eta})} \ \lesssim \ C_\eta\, 
	N^{\frac{1}{2} + \eta}\,
        \lp{u_{1,N}(0)}{L^2} \ , \label{unit_ang_str_est}
\end{equation}
where the implicit constants in the above inequality depend only on $n$ and $r$. 
\end{prop}\ret

\noindent In the case of $n=2$ spatial dimensions, we will also prove
the following:\\

\begin{prop}[Endpoint $(2+1)$ dimensional 
unit frequency Strichartz estimate for angularly regular
data]\label{n=2_unit_ang_str_th}
Let $n=2$ be the number of spatial dimensions, and let $u_{1,N}$ be a
unit frequency, angular frequency localized
solution to the homogeneous wave equation $\Box u_{1,N} = 0$.  Then
for every $0 < \eta$ there is a $2 < q_\eta$ such that $q_\eta \to 2$ as
$\eta\to 0$ and such that the following estimate holds:
\begin{equation}
        \lp{u_{1,N}}{L^{q_\eta}(L^\infty)} \ \lesssim \ C_\eta\, 
	N^{\frac{1}{2} + \eta}\,
        \lp{u_{1,N}(0)}{L^2} \ . \label{n=2_unit_ang_str_est}
\end{equation}
\end{prop}\ret

\noindent Assuming now the validity of Proposition \ref{unit_ang_str_th},
we prove Theorem \ref{full_ang_str_th} as follows: Given exponents $(q,r)$
such that $\frac{1}{q}+\frac{\sigma}{r} \geqslant \frac{\sigma}{2}$ and
$\frac{1}{q}+\frac{\sigma_\Omega}{r} < \frac{\sigma_\Omega}{2}$,
we first reduce things to the case where $q=2$. We define $t$ and
$r_1$ via the formulas:
\begin{align}
	\frac{1}{q} \ &= \ \frac{t}{2} \ , \notag\\
	\frac{1}{r} \ &= \ \frac{t}{r_1} + \frac{t-1}{2} \ . \notag
\end{align}
Notice that we can find such a $0\leqslant t \leqslant 1$ and
$\frac{2(n-1)}{n-2} < r_1$ due to the range of $(q,r)$. Therefore,
interpolating with the energy estimate (using Proposition 
\ref{L2_Omega_interp_prop}), it suffices to prove 
\eqref{full_ang_str_est} for indices $(2,r_1)$ with $\frac{2(n-1)}{n-2} < r_1$.\\

Next, using the decomposition \eqref{ang_lp_sp_sum} an rescaling the
spatial frequency of each term in the resulting sum,
it suffices to show that:
\begin{equation}
	\lp{u_{1,N}}{L^2(L^{r_1})} \ \lesssim \ C_\epsilon\,
	N^{(1+\epsilon)(\frac{n-1}{r_1} - \frac{n-3}{2})}\,
	\lp{u_{1,N}(0)}{L^2} \ . \label{red_ang_str_est}
\end{equation}
To do this, we choose:
\begin{equation}
	0\leqslant \frac{2(n-1)}{r_1}-(n-3) < t_\epsilon < 1 \ , \label{t_epsl_choice}
\end{equation}
such that there exists an $0 < \eta$ with the property that:
\begin{equation}
	(\frac{1}{2} + \eta)t_\epsilon \ \leqslant \
	(1+\epsilon)(\frac{n-1}{r_1} - \frac{n-3}{2}) \ . \label{desired_N_bnd}
\end{equation}
That such choices are possible follows from our assumptions on the range
of $r_1$ and the identity:
\begin{equation}
	(\frac{1}{2}+\frac{\epsilon}{2})t_0  \ = \  (1+\epsilon)(\frac{n-1}{r_1}
	- \frac{n-3}{2}) \ , \notag
\end{equation} 
where we have set $t_0 = \frac{2(n-1)}{r_1} - \frac{n-3}{2}$. 
Because of the range \eqref{t_epsl_choice}, we see that is is also possible
to choose an $\frac{2(n-1)}{n-2} < r_{\eta_*}\leqslant \frac{2(n-1)}{n-3}$
with the property that:
\begin{equation}
	\frac{t_\epsilon}{r_{\eta_*}} + \frac{1-t_\epsilon}{2(n-1)/(n-3)}
	\ = \ \frac{1}{r_1} \ . \notag
\end{equation}
Furthermore, using Proposition \ref{unit_ang_str_th} and possibly a Sobolev
embedding, we see that we have the estimate:
\begin{equation}
	\lp{u_{1,N}}{L^2(L^{r_{\eta_*}})} \ \lesssim \ C_{\eta, \eta_*}\,
	N^{\frac{1}{2} + \eta}\, \lp{u_{1,N}(0)}{L^2} \ . \notag
\end{equation}
Interpolating this last line with the $L^2(L^\frac{2(n-1)}{n-3})$ endpoint
of \eqref{str_est} we have achieved \eqref{red_ang_str_est}. Therefore, in
the sequel, we will concentrate on the proof of \eqref{unit_ang_str_est}.\\ \\


\subsection{The Hankel transform}
Finally, to wrap things up for this section, we record here the following
formula for the action of the inverse Fourier transform on the decomposition:
\begin{equation}
        \widehat{f} \ = \ \sum_{l,i} \ \widehat{c^l_i} \, Y^l_i \ . \notag
\end{equation} 
As is well known, this is given by a series of Hankel transforms. The formula
is (see \cite{SWeuclid}):
\begin{equation}
        f(r\omega) \ = \ \sum_{l,i} \ 2\pi(\sqrt{-1})^{l} \ r^\frac{2-n}{2} \ 
        \int_0^\infty  \ J_{\frac{n-2}{2} + l }\, (2\pi r \rho) \ 
        \widehat{c^l_i}(\rho)
        \ \rho^\frac{n}{2}d\rho \ \cdot \ Y^l_i(\omega) \ . \label{hankel_exp}
\end{equation}
Here $J_{s}(y)$ is the Bessel function of order $s$. For $-\frac{1}{2} < s$, 
this is given by the integral formula:
\begin{equation}
        J_s(y) \ = \ \frac{(y/2)^s}{\Gamma[(2s+1)/2]\, \Gamma(1/2)} \ 
        \int_{-1}^1 \ e^{ity}\, (1 - t^2)^\frac{2s-1}{2} \ dt \ . 
        \label{bessel_s}
\end{equation}

\ret \ret

\section{The Hankel--$\varphi$ transform}

As we see from the formula \eqref{hankel_exp} of the last subsection,
it is possible to expand the expression\footnote{Throughout this subsection,
we will work with the operator $e^{-i t \sqrt{-\Delta}}$ instead of
$e^{i t \sqrt{-\Delta}}$. Of course this is just a matter of notational 
convenience, as can be seen for instance by time reversal.}
$e^{-i t \sqrt{-\Delta}}f$
in terms of spherical harmonics as:
\begin{equation}
        e^{ - i t \sqrt{-\Delta}}f\, (r\omega) \ = \ 
        \sum_{l,i} \ c^l_i(t,r)\, Y^l_i(\omega) \ , \label{f_wv_prop_exp}
\end{equation}
where the coefficients $c^l_i$ are given by the Hankel transform formula:
\begin{equation}
        c^l_i(t,r) \ = \ 2\pi(\sqrt{-1})^{l} \ r^\frac{2-n}{2} \ 
        \int_0^\infty  \ J_{\frac{n-2}{2} + l }\, (2\pi r \rho) \ 
        e^{ - 2\pi i t \rho}\ \widehat{c^l_i}(\rho)
        \ \rho^\frac{n}{2}d\rho \ . \label{Hankel_wave_coef}
\end{equation}
Here, as in the previous subsection, the $\widehat{c^l_i}(\rho)$ are the 
coefficients in the spherical harmonic expansion of $\widehat{f}\, (\rho)$.
Also, the coefficients $c^l_i(t,r)$ should not be confused with the inverse
Fourier transform of the $\widehat{c^l_i}(\rho)$.\\

We would now like to be able to localize the expressions 
\eqref{Hankel_wave_coef}
in a manner analogous to the localization of the integral 
\eqref{sph_wave_osc_int}. This would be a relatively simple matter, if we could
show that the asymptotic \eqref{bessel_asym} held uniformly in $n$. That
is, if there was an asymptotic of the form \eqref{bessel_asym} for
$J_{\frac{n-2}{2} + l }\, (y)$ which held uniformly as $l\to \infty$. 
Unfortunately, it is well known that 
this is only the case for the region $l \ll \sqrt{y}$ (see 
\cite{Wbessel}). In the transition regions, that is when $ \sqrt{y} \lesssim 
l \lesssim y$, the asymptotic for $J_{\frac{n-2}{2} + l }\, (y)$ becomes
quite complicated. Roughly speaking, it begins to loose oscillations
in $y$ while it gains decay in the parameter $l$. Because of this, it does
not seem feasible to try and compute an approximate formula for
$J_{\frac{n-2}{2} + l }$ and then substitute it into the integrals
\eqref{Hankel_wave_coef}. Instead we will use a more straight forward
approach, by first localizing the $\widehat{c^l_i}$ in frequency as a Fourier 
series, 
just as we had done for $\widehat{f_1}$ in the integral 
\eqref{sph_wave_osc_int},
and then computing the integral \eqref{Hankel_wave_coef} directly by using  
appropriate integral representations for the $J_{\frac{n-2}{2} + l }$.\\

Since we are assuming that the initial data in Theorem \ref{unit_ang_str_th}
is unit frequency, we will assume that all of the coefficient functions
$\widehat{c^l_i}(\rho)$ in the integrals \eqref{Hankel_wave_coef} are
supported on the interval $(\frac{1}{2} , 2)$. We may take their Fourier
series developments on the interval $(0,4)$, and we record these as:
\begin{equation}
        \widehat{c^l_i}(\rho) \ = \ \sum_k \ c^l_{i,k} \, e^{i \frac{\pi}{2}
        \rho} \ . \label{c_fourier_series}
\end{equation}
Expanding the integral \eqref{Hankel_wave_coef} in terms of the above formula,
we see that:
\begin{equation}
        c^l_i(t,r) \ = \ 2\pi(\sqrt{-1})^{l} \sum_k \ r^\frac{2-n}{2} \
        c^l_{i,k} \ \psi^l_{  t  - \frac{k}{4} } (r) \ , \label{c_wavelet_devel}
\end{equation}
where:
\begin{equation}
        \psi^l_{ t - \frac{k}{4} } (r) \ = \
        \int_0^\infty  \ J_{\frac{n-2}{2} + l }\, (2\pi r \rho) \ 
        e^{-2\pi i ( t - \frac{k}{4})\rho }\ \chi_{(\frac{1}{4} , 4)}(\rho)
        \ d\rho \ . \label{hankel_wavelet}
\end{equation}
In the above formula $\chi_{(\frac{1}{4} , 4)}$ is a smooth bump function
on the interval $(\frac{1}{4} , 4)$. Notice that this is not necessarily
equal to $1$ on any interval because we have absorbed the volume element
into our definition of $\chi_{(\frac{1}{4} , 4)}$.
We call the right hand side of \eqref{c_wavelet_devel} the 
\emph{Hankel--$\varphi$} transform of the function $c^l_i(t,r)$. 
We would now like to be able to give a precise bound on how well localized
the functions $\psi^l_m$ are in physical space
for the various values of the half-integer
parameter $l$ and the real variable $m$. This brings us to the main result of 
this subsection:\\

\begin{prop}[Asymptotics of the functions $\psi^l_{m}$]\label{psi_asym_prop}
Let $\psi^l_{m}$ be the function given by the formula \eqref{hankel_wavelet}
for $m = t - \frac{k}{4}$. Then for every set integers $0 \leqslant N_1 , N_2$, 
there exists a
constant $C_{N_1,N_2}$ depending only on the $N_i$ (and the dimension $n$) 
such that 
the following asymptotics hold uniform in the parameters $l$, $m$, and $r$:
\begin{align}
        |\psi^l_{m} (r) | \ &\leqslant  \ r^{\frac{n-2}{2}} 
	\cdot \frac{  C_{N_1,N_2}  }{(1 + |m|)^{N_1}}\cdot\left(\frac{1}{1 + l}
	\right)^{N_2}\ ,&
        0 \ \leqslant \ r \ < \ 1 \ &, \label{psi_asym0} \\
        |\psi^l_{m} (r) | \ &\leqslant \
        \frac{C_{N_1,N_2}}{(1 + r + |m|)^\frac{1}{2}
	\cdot (1 + \big| r  - |m| \big|)^{N_1}}
	\cdot\left(\frac{r^\frac{1}{2}}{1 + l}\right)^{N_2}
        \ ,&  1 \ \leqslant \ r \ \leqslant \ |m|+1 &. \label{psi_asym1}\\
	|\psi^l_{m} (r) | \ &\leqslant \ 
        \frac{1}{ ( r^2 - m^2 )^\frac{1}{2}}\cdot R(l,m,r)\ ,&
        |m| + 1 \ < \ r \ &. \label{psi_asym3}
\end{align}
The extra term $R(m,r)$ in line \eqref{psi_asym3} above is a positive 
function with the bound:
\begin{equation}
        R(l,m,r) \ \leqslant \
        C_{N_1,N_2}\left(\
        \frac{1}{ (1 + \big| |m| - r \big|)^{N_1}} \ + \ 
        \min_\pm\Big\{ 
        \Big( \frac{l}{ \sqrt{r^2 - m^2}} \Big)^{\pm N_2}  \Big\}\ \right) \ . 
        \label{psi_asym3_remainder}
\end{equation}
\end{prop}\ret

\begin{rem}\label{coh_rem}
The downside of the above asymptotic is of course the region governed by
\eqref{psi_asym3}. When $l^2 \lesssim r$, one can see that the extra factor
\eqref{psi_asym3_remainder} will allow this asymptotic to look like \eqref{psi_asym1}.
Notice that this is consistent with the fact that one has the asymptotic
\eqref {bessel_asym} for $J_s(y)$ in this region, and is what is responsible for
the good localization \eqref{psi_asym} for spherically symmetric waves. Unfortunately,
it does not seem like one can do much to improve \eqref{psi_asym3} in the region
where $r\ll l^2$ (except for the extra factor of $R(l,m,r)$). In fact,
if one assumes there is an asymptotic for $\psi^l_{m} (r)$ in this region
which is of the form \eqref{psi_asym1},
and one sets the $c^l_{i,k}\equiv 1$ in the sum \eqref{c_wavelet_devel}
for a fixed $l,i$, by putting absolute values around the sum \eqref{c_wavelet_devel}
one would get an asymptotic that looks like: 
$|c^l_i(0,r)| \lesssim \frac{1}{r^\frac{1}{2}}$. But in this case, $c^l_i(t,r)$
corresponds to a delta function along the radial variable in Fourier space,
say supported at the point $\frac{1}{2\pi}$. Therefore we would have shown
a bound like $|J_s(r)| \lesssim \frac{1}{r^\frac{1}{2}}$ uniform in $s$. This
violates the well known asymptotic for Bessel functions: $|J_s(s)|\sim s^{-\frac{1}{3}}$
(see \cite{Wbessel}). It would be interesting to know if there is a more coherent
decomposition of the Hankel transform that could eliminate this problem. 
\end{rem}\ret


\begin{proof}[proof of Proposition \ref{psi_asym_prop}]
The asymptotics \eqref{psi_asym0}--\eqref{psi_asym1} follow
more or less directly from appropriate integral formulas for the
$J_s(y)$. We will need to split the proof into the two cases: \ \{$r\leqslant 1$ or 
$1 < r \leqslant |m| + 1$\} and
\{$|m|+1 < r $\}.\\ \\

\noindent \textbf{Case 1: \ $r \leqslant 1$ or $1 < r \leqslant |m| + 1$.}\\
Here we use a standard integral representation for Bessel functions which differs
from \eqref{bessel_s}.
For $s\in\mathbb{N}$, one has the following formula:
\begin{equation}
        J_s(y) \ = \ \frac{(-i)^s}{2\pi}\ \int_0^{2\pi} \ e^{i y \cos \theta}
        e^{-is \theta} \ d\theta \ .  \label{J_s_int_form}
\end{equation}
This can be proved by a simple recursive argument (see \cite{SWeuclid},
Chapter 4, Lemma 3.1). All one has to is to show that for both the 
integral formulas \eqref{J_s_int_form} and \eqref{bessel_s},
the function $J_s(y)$ satisfies the recursive relation:
\begin{align}
        \frac{d}{dt} \left[ t^{-s} J_s(y) \right] \ = \ t^{-s}
        J_{s + 1}(y) \ ,& &0 < t \ , \label{J_rec_form}
\end{align}
for $s\in \mathbb{N}$. In light of \eqref{J_rec_form}, the equality of 
\eqref{J_s_int_form} and \eqref{bessel_s} is reduced to
showing that it is true when $s=0$. This can be achieved directly
through a change of variables.\\

Now using periodicity, integrating over an adjacent interval
of length $2\pi$, and averaging, we see that for $s\in \mathbb{N}$
the following integral representation also holds:
\begin{equation}
        J_s(y) \ = \ 
        \frac{(-i)^s}{4\pi} \ \int_{-2\pi}^{2\pi} \ e^{i y \cos \theta}
        e^{-is \theta} \ d\theta . \label{J_s_half_int_form}
\end{equation}
Moreover, by a direct calculation, its not hard to see that the recursive
relation \eqref{J_rec_form} is satisfied by both the integrals 
\eqref{J_s_half_int_form} and \eqref{bessel_s} whenever
$s \in \frac{1}{2}\cdot \mathbb{N}$. Therefore, throughout the sequel,
we may assume that \eqref{J_s_half_int_form} is our definition of the
Bessel function that appears in the integral formula \eqref{hankel_wavelet}
for the function $\psi^l_m$ in dimension $n$. Making this substitution
yields:
\begin{align}
        \psi^l_m (r)  \ &= \ \frac{(-i)^s}{4\pi} \ 
        \int_{-\infty}^\infty \int_{-2\pi}^{2\pi} \ 
        e^{2\pi i (r \cos\theta - m)\rho }\
        e^{-i(\frac{n-2}{2} + l) \theta} \ \chi_{(\frac{1}{4} , 4)}(\rho)
        \ d\theta\, d\rho \ , \notag \\
        &= \ \frac{(-i)^s}{4\pi} \ \int_{-2\pi}^{2\pi} \ 
        \check\chi_{(\frac{1}{4} , 4)} (r \cos\theta - m) \ 
        e^{-i(\frac{n-2}{2} + l) \theta}\ d\theta \ . \label{working_psi_form}
\end{align}\ret

We begin by proving the asymptotic \eqref{psi_asym0}. In fact, we will prove a bit
more. We will show that the asymptotic \eqref{psi_asym0} holds for
$r\leqslant 30$. Our first step
will be to pick up the decay in the $(1+l)$ parameter by integrating 
by parts the expression \eqref{working_psi_form} $N_2$ times. The resulting
expression looks like:
\begin{equation}
       \psi^l_m (r) \ = \ \left(\frac{i}{\frac{n-2}{2} + l}\right)^{N_2}
       \ \sum_{k=1}^{N_2}\ \int_{-2\pi}^{2\pi} \ r^k p_k(\theta)\
       \check\chi^{(k)}_{(\frac{1}{4} , 4)} (r \cos\theta - m) \ 
        e^{-i(\frac{n-2}{2} + l) \theta}\ d\theta \ , \label{int_parts_exp}
\end{equation}
where the $p_k(\theta)$ in the above formula denote some specific trigonometric
polynomials of degree $k$ who's exact form is not important for our analysis. 
The next step is to gain the decay in $r$ in
conjunction with the damping in terms of inverse powers of $(1+|m|)$. To get
this, we Taylor expand each 
$\check\chi^{(k)}_{(\frac{1}{4} , 4)} (h(\theta) - m)$, where
\begin{equation}
        h(\theta) \ := \ r \cos\theta  \ , \label{h_def}
\end{equation}
around the point $h=0$. Notice that this is consistent with the 
fact that we are investigating
the region where $r$ is bounded. We now define the dimensional constant 
$M =\lceil\frac{n-2}{2}\rceil$, 
we record this Taylor expansion as:
\begin{equation}
        \check\chi^{(k)}_{(\frac{1}{4} , 4)} (h - m) \ = \ 
	\sum_{j=0}^{M-1} \
	\frac{1}{j!}\, \check\chi^{(k+j)}_{(\frac{1}{4} , 4)} (- m)\cdot
	h^j \ + \ \frac{1}{M!}\, \check\chi^{(k+M)}_{(\frac{1}{4} , 4)} 
	(u(h) - m)\cdot h^{M} \ , \label{taylor_exp}
\end{equation}
where $u(h)$ is some smooth function such that $|u(h)|\leqslant 30$.
Substituting the Taylor expansion \eqref{taylor_exp} into the integral
\eqref{int_parts_exp}, we see that we may write:
\begin{equation}
        \psi^l_m (r) \ = \ A + B \ , \notag
\end{equation}
where:
\begin{equation}
        A \ = \ \left(\frac{i}{\frac{n-2}{2} + l}\right)^{N_2}\
	\sum_{k + j \ < \ M}\
	\frac{1}{j!}\, \check\chi^{(k+j)}_{(\frac{1}{4} , 4)} (- m)\
	\int_{-2\pi}^{2\pi} \
	r^{j+k} p_k(\theta)\cos^j(\theta)\ 
	e^{-i(\frac{n-2}{2} + l) \theta}\ d\theta \ . \notag
\end{equation}
and: 
\begin{equation}
        |B| \ \leqslant \ C_{N_2}\cdot r^M\cdot
	\left(\frac{1}{\frac{n-2}{2} + l}\right)^{N_2}
	\cdot\ \sum_{k=M}^{N_2 +M}\ 
	\sup_{m-30 \leqslant x \leqslant m+30}\, |
	\check\chi^{(k)}_{(\frac{1}{4} , 4)}(x)| \ . \notag
\end{equation}
We can further estimate the term $B$ above by using the fact that 
$\check\chi$ and all of its derivatives have rapid
decay away from the origin. This follows immediately from 
the fact that $\check\chi$ is the inverse Fourier transform
of a smooth $O(1)$ bump function. We record this observation as:
\begin{equation}
        | \check\chi^{(k)}_{(\frac{1}{4} , 4)} (y) | \ \lesssim \
        \frac{C_{k,N_1}}{( 1 + |y|)^{N_1}} \ . \label{check_chi_bound}
\end{equation}
By adding things up and introducing a large enough constant, this
allows us to write:
\begin{equation}
        |B| \ \leqslant \ C_{N_1,N_2}\cdot r^M\cdot
	\left(\frac{1}{\frac{n-2}{2} + l}\right)^{N_2}
	\cdot \frac{1}{(1+|m|)^{N_1}} \ . \notag
\end{equation}
Recalling now that we have set $M=\lceil\frac{n-2}{2}\rceil$
and that we also have $3\leqslant n$,
we see that in order to achieve the bound \eqref{psi_asym0},
all we need to do is to control the expression for $A$. This is easy
to do because a moments inspection shows that in fact one has $A\equiv 0$.
This can be readily seen for even dimensions, that is when $n$ is even,
because in this case the trigonometric polynomials
under the integral sign in the expression for $A$ are of degree
strictly less than $\frac{n-2}{2}$. By orthogonality,
the whole expression then integrates to zero. In the case of odd dimension,
each integral is still zero thanks to the fact that $\frac{n-2}{2}$
is a half integer expression, where as the term:
\begin{equation}
        p_k(\theta)\cos^j(\theta)\ e^{-il} \ , \notag
\end{equation}
is a trigonometric polynomial of integer degree. Since we are integrating
over the double torus, $[-2\pi,2\pi]$, a rescaling turns the expression under the integral
sign in $A$ into a product of
even degree trigonometric polynomials and odd degree trigonometric
polynomials. Therefore one has the needed orthogonality. This completes the proof
of \eqref{psi_asym0}.\\ \\

We now turn our attention to proving the asymptotic
\eqref{psi_asym1} for the regime where $1 < r \leqslant |m|+1$. 
For the remainder of this section we will assume that $m$ is positive, as the other
case can be dealt with by a similar argument. Furthermore,
using the bound that we proved in the previous discussion, we can without loss of
generality assume that $20 \leqslant m$.\\

Our first step is to split the integral on the right hand side of \eqref{working_psi_form}
smoothly into the regions where $(1-\cos\theta) \ll 1$ and otherwise. 
To realize this split, we restrict the integral \eqref{working_psi_form} to the regions:
\begin{align}
        R_1 \ &= \ \{\theta \big| |\theta| < 1\} \ , \notag \\
	R_2 \ &= \ \{\theta\big| |\theta\pm2\pi| < 1\} \ , \notag \\
	R_3 \ &= \ [-2\pi,2\pi]\setminus (\{\theta \big| |\theta| < \frac{1}{2}\}\cup
	\{\theta\big| |\theta\pm2\pi| < \frac{1}{2}\}) \ . \notag
\end{align}
Notice that by symmetry, we only need to consider the regions $R_1$ and $R_3$. On these 
regions, a bit of explicit computation using Taylor expansions shows that:
\begin{align}
        h(\theta) - m \ &= \ -\frac{1}{2}r [u(\theta)]^2 - (m-r) \ ,&
	&\theta\in R_1 \ , \label{h_bound_in_R1} \\
	|h(\theta) - m| \ &\geqslant \ \frac{1}{100}m \ ,&
	&\theta\in R_2 \ \hbox{and} \ , \label{h_bound_in_R2}
\end{align}
where the function $u(\theta)$ in the first line above satisfies the bound:
\begin{equation}
        \frac{1}{2} \ < \ u(\theta)/ \theta^2 \ < \ 2 \ . \label{bound_for_u}
\end{equation}
Notice that to get \eqref{h_bound_in_R2} we have used the condition 
$20 \leqslant m$.
Next, using another simple calculation involving the Taylor
series of trigonometric functions, as well as the estimates \eqref{check_chi_bound},
\eqref{h_bound_in_R1} and \eqref{bound_for_u}, we see that one may write
for $|\theta| < 1$:
\begin{align}
        \Big| \ \frac{d^{N_2}}{d\theta^{N_2}}\left[
	\check\chi(h(\theta)-m)\right] \ \Big| \ &\leqslant \ C_{N_2}
	\ \sup_{0\leqslant j\leqslant N_2} |\check\chi^{(j)}
	(h(\theta)-m)|\cdot 
	\sum_{k=0}^{ \lfloor\frac{N_2}{2}\rfloor} \ r^{N_2-k}|\theta|^{N_2-2k}
	\ , \notag \\
	&\leqslant \ \frac{C_{N_1,N_2}}{( 2 + m - r + r\theta^2)^{N_1 + N_2 + 1}}
	\cdot \sum_{k=0}^{\lfloor\frac{N_2}{2}\rfloor} 
	\ r^{N_2-k}|\theta|^{N_2-2k} \ .
	\label{chi_derivative_bound}
\end{align}
We are now in a position to bound the integral 
\eqref{working_psi_form} on the region $R_1$.
To do this, we will employ a smooth cutoff function $\chi_{(-1,1)}$
such that:
\begin{equation}
        \chi_{(-1,1)}(\theta) \ = 
	\begin{cases}
	        1 \ ,& \theta\in (-\frac{1}{2},\frac{1}{2}) \ , \notag \\
		0 \ ,& \theta\notin (-1,1) \ . \notag
	 \end{cases} \notag
\end{equation}
We are trying to prove the estimate \eqref{psi_asym1} for the expression:
\begin{equation}
       \int \  \chi_{(-1,1)}(\theta)\cdot
        \check\chi_{(\frac{1}{4} , 4)} ( h(\theta)-m ) \ 
        e^{-i(\frac{n-2}{2} + l) \theta}\ d\theta \ . \notag
\end{equation}
Integrating by parts $N_2$ times in the above integral and
using the bound \eqref{chi_derivative_bound} (along with the fact that we are in dimension 
$3\leqslant n$) we compute:
\begin{align}
        &\Big| \int \  \chi_{(-1,1)}(\theta)\cdot
        \check\chi_{(\frac{1}{4} , 4)} (h(\theta)-m) \ 
        e^{-i(\frac{n-2}{2} + l) \theta}\ d\theta \Big| \ , \notag \\
	\leqslant \ &\left(\frac{1}{l+1}\right)^{N_2}\cdot
	\Big| \int \  \frac{d^{N_2}}{d\theta^{N_2}}\left[
	\chi_{(-1,1)}(\theta)\cdot
        \check\chi_{(\frac{1}{4} , 4)} (h(\theta)-m)\right] \ 
        e^{-i(\frac{n-2}{2} + l) \theta}\ d\theta \Big| \ , \notag \\
	\leqslant \ 
	&\left(\frac{1}{l+1}\right)^{N_2}\cdot\sum_{0\leqslant k \leqslant \frac{N_2}{2}}
	\ \int \ \frac{ C_{N_1,N_2}}{( 2 + m - r + r\theta^2)^{N_1 + N_2 + 2}}
	\cdot r^{N_2-k}|\theta|^{N_2-2k} d\theta \ , \notag \\
	\leqslant \ &C_{N_1,N_2} \ \frac{1}{r^\frac{1}{2}}\cdot
	\left(\frac{r^\frac{1}{2}}{l+1}\right)^{N_2}
	\ \int \ \frac{ 1 + |\theta|^{N_2} }{( 2 + m - r + 
	\theta^2)^{N_1 + N_2 + 2}}
	\ d\theta \ , \notag \\
	\leqslant &\ \frac{C_{N_1,N_2}}
	{r^\frac{1}{2}( 2 + m -r )^{N_1+1}}\cdot
	\left(\frac{r^\frac{1}{2}}{l+1}\right)^{N_2}
	\ \int \ \frac{ 1 + |\theta|^{N_2} }{( 1 + \theta^2)^{N_2 + 1}}
	\ d\theta \ , \notag \\
	\leqslant &\ \frac{C_{N_1,N_2}}
	{( 1 + m + r)^\frac{1}{2}( 2 + m -r )^{N_1}}\cdot
	\left(\frac{r^\frac{1}{2}}{l+1}\right)^{N_2} \ . \notag
\end{align}
Thus, we have proved \eqref{psi_asym1} for this portion of things.\\ 

It remains to prove the estimate \eqref{psi_asym1} for the 
region $R_3$. We suggestively (c.f. \eqref{h_bound_in_R2})
denote the cutoff here by $\chi_{m \lesssim |h(\theta)| }= 
(1- \chi_{(-1,1)})$. The calculation is
essentially the same as what was done above, except that here we can afford
to be more careless about the powers of $r$ which come up through
integration by parts. Integrating by parts $N_2$ times and using the bound 
\eqref{h_bound_in_R2} in conjunction with the estimate \eqref{check_chi_bound},
we see that:
\begin{align}
        &\Big| \int \  \chi_{ m \lesssim |h(\theta)|  }(\theta)\cdot
        \check\chi_{(\frac{1}{4} , 4)} (h(\theta)-m) \ 
        e^{-i(\frac{n-2}{2} + l) \theta}\ d\theta \Big| \ , \notag \\
	\leqslant \ &C_{N_2}\ \left(\frac{1}{l+1}\right)^{N_2}\cdot
	\sum_{0\leqslant k \leqslant N_2} \ 
	r^{N_2}\cdot \sup_{\frac{1}{100}m \leqslant |h - m|}
	|\check\chi^{(k)}_{(\frac{1}{4} , 4)} (h-m)|
	\ , \notag \\
	\leqslant \ &C_{N_1,N_2}\ \left(\frac{1}{l+1}\right)^{N_2}\cdot
	\frac{r^{N_2}}{(1 + m )^{N_1 + \frac{N_2 + 1}{2}}} \ , \notag \\
	\leqslant &\ \frac{C_{N_1,N_2}}
	{( 1 + m + r)^\frac{1}{2}( 1 + m -r )^{N_1}}\cdot
	\left(\frac{r^\frac{1}{2}}{l+1}\right)^{N_2} \ . \notag
\end{align}
This completes the proof of \eqref{psi_asym1} and ends the
demonstration of \textbf{case 1}.\\ \\

\noindent \textbf{Case 2: \ $|m| + 1< r$.}\\
In this section, the main difficulty will be for us to 
incorporate the remainder term 
\eqref{psi_asym3_remainder} into our asymptotic. In order to motivate the
steps we will take here, we argue heuristically as follows: 
We first define the auxiliary function:
\begin{equation}
	g(\theta) \ := \ r\cos\theta - m \ . \notag
\end{equation}
It is clear that the main contribution
to the integral \eqref{working_psi_form} integrating around points where the 
function $g(\theta)$ defined above vanishes. By symmetry, we need only
consider such points which are positive. We'll call the one closest to
zero $\theta_0$. Now, expanding $g(\theta)$ around this zero gives:
\begin{equation}
        h(\theta) \ = \ -\sqrt{r^2 - m^2}(\theta-\theta_0) \ - \ \frac{m}{2}\, 
        (\theta - \theta_0)^2 \ + \ O\big( (\theta -\theta_0)^3\big) 
        \ , \label{h0_taylor_exp}
\end{equation}
where again $0 < \theta_0 = \cos^{-1} (\frac{m}{r})$. Therefore, it is natural
to expect that:
\begin{align}
        \hbox{(R.H.S.)}\eqref{working_psi_form} \ &\approx \ 
        \frac{(-i)^s}{4\pi} \ \int_{-2\pi}^{2\pi} \ 
        \check\chi_{(\frac{1}{4} , 4)} \big(  
        -\sqrt{r^2 - m^2}(\theta-\theta_0)   \big) \ 
        e^{-i(\frac{n-2}{2} + l) \theta}\ d\theta \ , \notag \\
        &\approx \ \frac{e^{-i(\frac{n-2}{2} + l)\theta_0}  }
        {(r^2 - m^2)^\frac{1}{2}} \ \int_{-\infty}^{\infty} \ 
        \check\chi_{(\frac{1}{4} , 4)}(\theta) \ 
        e^{i\, \frac{(n-2)/2 + l}{\sqrt{r^2-m^2}} \cdot \theta} \ d\theta
        \ , \notag \\
        &= \ \frac{e^{-i(\frac{n-2}{2} + l)\theta_0}  }{(r^2 - m^2)^\frac{1}{2}}
        \cdot \left(i\, \frac{(n-2)/2 + l}{\sqrt{r^2-m^2}}\right)^{N_2}
         \ \int_{-\infty}^{\infty} \ 
        \check\chi_{(\frac{1}{4} , 4)}^{(-N_2)}(\theta) \ 
        e^{i\, \frac{(n-2)/2 + l}{\sqrt{r^2-m^2}} \cdot \theta} \ d\theta
        \ , \label{approx_working_psi}
\end{align}
where the $\check\chi_{(\frac{1}{4} , 4)}^{(-N_2)}$, 
$N_2= 0 , \pm 1 , \pm 2 , \ldots$, 
are the (unique) set of derivatives and derivatives of 
$\check\chi_{(\frac{1}{4} , 4)}$ which vanish at infinity. These 
satisfy the bounds similar to \eqref{check_chi_bound}, even if $N_2$ is positive:
\begin{equation}
        \check\chi_{(\frac{1}{4} , 4)}^{(-N_2)}(y) \ \leqslant \
        \frac{C_{N_2}}{(1 + |y|)^2} \ . \label{chi_deriv_bnd}
\end{equation}
This, of course, comes from the fact that $\check\chi_{(\frac{1}{4} , 4)}$
is the inverse Fourier transform of a unit frequency function, and is therefore
supported away from the origin in Fourier space.
Using  \eqref{chi_deriv_bnd} we easily get the following bound for arbitrary
non-negative $N_2$:
\begin{equation}
        \big| \hbox{(R.H.S.)}\eqref{approx_working_psi} \big|
        \ \leqslant \ \frac{C_{N_2}}{(r^2 - m^2)^\frac{1}{2}}\cdot
        \min_\pm \left\{ \left(\frac{l}{\sqrt{r^2-m^2}}\right)^{\pm N_2}
        \right\} \ . \notag
\end{equation}\\

In order to make the previous argument rigorous, we need to justify
the two approximations used in the lines directly above 
\eqref{approx_working_psi}. We will not be able to do this 
completely, which will be responsible for the extra term in the 
formula \eqref{psi_asym3_remainder} for $R(l,m,r)$.
As in the previous two subsections, we argue by isolating  the
interval of integration. By symmetry, and the decay bound 
\eqref{check_chi_bound},
we can without loss of generality assume that we are integrating 
\eqref{working_psi_form} over the interval $[0,\frac{3\pi}{4}]$. We will
first need to go a little further and chop some more off of the left hand side
of this interval. What we'll do is take $\theta_1$  (we're still keeping
$\theta_0 = \cos^{-1} (\frac{m}{r})$) to be such that:
\begin{equation}
        \frac{r}{\sqrt{r^2 - m^2}}\sin (\theta_1) \ = \ \frac{1}{2} \ . 
	\label{theta1_def}
\end{equation}
Taylor expanding $\sin\theta$ gives us the bound:
\begin{equation}
        \theta_1 \ \lesssim \ \frac{(r - m)^\frac{1}{2}}{r^\frac{1}{2}}
        \ . \label{theta1_bnd}
\end{equation}
Furthermore, we also have the bound:
\begin{align}
        \frac{3}{4}(r-m) \ \leqslant \ \frac{\sqrt{3 r^2 + m^2}}{2} - m
        \ \leqslant \ r \cos\theta - m \ = \ g(\theta) \ ,&
        &\theta \ \in \ [0,\theta_1] \  . \ \label{short_interval_h_bnd}
\end{align}
Now, using \eqref{theta1_bnd} and \eqref{short_interval_h_bnd} 
above, in conjunction with the
asymptotic \eqref{check_chi_bound}, we immediately see that:
\begin{equation}
        \int_0^{\theta_1} \ \big| \check\chi_{(\frac{1}{4} , 4)}
        \big(g(\theta)\big) \big| \ d\theta \ \lesssim \
        \frac{(r - m)^\frac{1}{2}}{r^\frac{1}{2}}\cdot
        \frac{C_{N_1}}{(r - m)^{N_1}} \ . \notag
\end{equation}
This is enough to give \eqref{psi_asym3} for this portion of things.
Therefore, it remains to compute the integral:
\begin{equation}
        I_1 \ = \ \int_{\theta_1}^\frac{3\pi}{4} \ 
        \check\chi_{(\frac{1}{4} , 4)}
        \big(g(\theta)\big) \ e^{-i(\frac{n-2}{2} + l) \theta}\ d\theta
        \ . \notag
\end{equation}
Keeping the Taylor expansion \eqref{h0_taylor_exp} in mind,
we now make the following change of variable for $I_1$:
\begin{equation}
        \varphi(\theta) \ = \ g\left(\frac{1}{\sqrt{r^2-m^2}}\, \theta + 
        \theta_0\right) \ . \notag
\end{equation}
Using this, we can write:
\begin{equation}
         | I_1 | \ = \ \frac{1}{\sqrt{r^2-m^2}}\ 
        \Big| \ \int_R \ 
        \check\chi_{(\frac{1}{4} , 4)}\big( \varphi(\theta) \big)
        \ e^{-i\frac{(n-2)/2 + l}{\sqrt{r^2-m^2}} \ \theta} \
        d\theta \ \Big| \ , \label{final_int}
\end{equation}
Where $R$ is the interval:
\begin{equation}
	R \ = \ \left[\sqrt{r^2-m^2}(\theta_1-\theta_0),
	\sqrt{r^2-m^2}(\frac{3\pi}{4} - \theta_0)\right] \ . \notag
\end{equation}
The desired result will
now follow by integrating by parts as many times as necessary the 
integral \eqref{final_int}. However, some care needs to be taken
in order to control terms involving $\varphi'$, which can be as big
as $\frac{r}{\sqrt{r^2-m^2}}$. Also, one needs to know that the 
higher derivatives of $\varphi$ possess some decay in order to 
deal with terms where the derivatives fall on $\frac{1}{\varphi'}$
instead of the exponential factor. We address these issues now.
The first observation we make is that we may integrate the bound:
\begin{align}
        \frac{1}{2} \ \leqslant \ \frac{r}{\sqrt{r^2 - m^2}}\sin \theta
        \ ,& &\theta \ \in \ [\theta_1,\frac{3\pi}{4}] \ , 
        \label{phi_prime_bound}
\end{align}
using the fact that $\varphi(0)=0$ to get that:
\begin{equation}
        |\theta| \ \lesssim |\varphi(\theta)| \ . \label{varphi_var_bnd}
\end{equation}
over the range of integration. However, this is not enough to 
address the issue of a product of the form 
$(\varphi')^k\cdot \check\chi^{(j)}_{(\frac{1}{4} , 4)}$. 
To handle this, notice that one has the bound:
\begin{align}
	\frac{r}{\sqrt{r^2-m^2}}\sin\theta \ \leqslant \
	m -r\cos\theta + 1 \ ,& \ 
	&\theta\in [\theta_1 , \frac{3\pi}{4}] \ , \label{varphi_der_bnd}
\end{align}
as can be seen from the fact that this bound is true for 
$\theta = \theta_1$, and that upon differentiating both sides of
this expression one is reduced to showing (the inequality for the 
derivatives is trivial for $\theta\in[\frac{\pi}{2},\frac{3\pi}{4}]$):
\begin{align}
	\frac{1}{\sqrt{r^2-m^2}}\ \leqslant \
	\tan\theta \ ,& &\theta\in [\theta_1 , \frac{\pi}{2}] \ . \notag
\end{align}
This then follows immediately from the increasing nature of $\tan\theta$
on this interval and from the identity $\tan(\theta_1) = 
\frac{\sqrt{r^2-m^2}}{m}$. Now, combining the bounds \eqref{varphi_der_bnd},
\eqref{varphi_var_bnd} and \eqref{check_chi_bound}, 
we see that for any positive integers $k,j$, one has the
estimates:
\begin{equation}
	(\varphi')^k(\theta)\cdot \check\chi^{(j)}_{(\frac{1}{4} , 4)}
	\big( \varphi(\theta) \big) \ \lesssim \
	\frac{C_{j,k}}{(1 + |\theta|)^2} \ . \label{bound_to_int}
\end{equation}
Finally, we record the fact that the higher derivatives of $\varphi$
satisfy the following simple bounds, which can be verified through
direct calculation:
\begin{align}
        | \varphi^{(k)} | \ \lesssim \ \frac{1}{ (r - m)^{k-1}} \ ,&
        &M \ = \ 2, 3 , \ldots \ . \label{phi_der_bnds}
\end{align} \\

We are now ready to bound \eqref{final_int}. First, we integrate as many
times as necessary, letting the derivatives fall on the term
$\check\chi_{(\frac{1}{4} , 4)}\big( \varphi(\theta) \big)$. The
resulting expression looks like:
\begin{align}
	&\Big| \ \int_R \ 
        \check\chi_{(\frac{1}{4} , 4)}\big( \varphi(\theta) \big)
        \ e^{-i\frac{(n-2)/2 + l}{\sqrt{r^2-m^2}} \ \theta} \
        d\theta \ \Big| \ , \notag \\
	\lesssim \ &\left(\frac{\sqrt{r^2-m^2}}{1 + l}\right)^{N_2}
	\int_R \ \Big| \frac{d^{N_2}}{d\theta^{N_2}}
	 \check\chi_{(\frac{1}{4} , 4)}\big( \varphi(\theta) \big)
	\Big| \ d\theta \ + \ 
	\sum_{k=0}^{N_2-1}\ |\frac{d^k}{d\theta^k}
	\check\chi_{(\frac{1}{4} , 4)}\big( \varphi(\theta) \big) |
	\ \Bigg|_{\partial R}\ . \notag
\end{align}
Expanding out the derivatives in the first term of the above
expression, and using the bounds \eqref{varphi_der_bnd} and 
\eqref{varphi_var_bnd}, we see that we can bound:
\begin{equation}
	\Big| \frac{d^{N_2}}{d\theta^{N_2}}
	 \check\chi_{(\frac{1}{4} , 4)}\big( \varphi(\theta) \big)
	\Big| \ \lesssim \ \frac{C_{N_2}}{(1 + |\theta|)^2}
	\ . \notag
\end{equation} 
Therefore, the integral in the first term above does not
cause us any trouble. As far as the boundary values are
concerned, we can use the bounds \eqref{theta1_def},
\eqref{short_interval_h_bnd}, \eqref{varphi_der_bnd}, and
\eqref{check_chi_bound} to show that:
\begin{align}
	&\sum_{k=0}^{N_2-1}\ |\frac{d^k}{d\theta^k}
	\check\chi_{(\frac{1}{4} , 4)}\big( \varphi(\theta) \big) |
	\ \Bigg|_{\partial R} \ , \notag \\ 
	\lesssim \
	&C_{N_2}\sup_{k_i\leqslant N_2}\left( 
	|\varphi^{(k_1)}(\theta)|^{k_2}
	\cdot|\check\chi^{(k_3)}_{(\frac{1}{4} , 4)}\big( 
	\varphi(\theta) \big)| \right)
	\Bigg|_{\sqrt{r^2-m^2}(\theta_1-\theta_0)}^
	{\sqrt{r^2-m^2}(\frac{3\pi}{4} - \theta_0)}
	\ , \notag \\
	\lesssim \ &C_{N_1,N_2}\left(\frac{1}{(r-m)^{N_1}} + 
	\frac{1}{(r+m)^{N_1}} \right) \ , \notag \\
	\lesssim \ &\frac{C_{N_1,N_2}}{(r-m)^{N_1}} \ . \notag
\end{align}\\

It remains to bound \eqref{final_int} where we integrate by
parts and let the derivatives fall on the exponential factor.
Doing this $N_2$ times yields:
\begin{align}
	&\Big| \ \int_R \ 
        \check\chi_{(\frac{1}{4} , 4)}\big( \varphi(\theta) \big)
        \ e^{-i\frac{(n-2)/2 + l}{\sqrt{r^2-m^2}} \ \theta} \
        d\theta \ \Big| \ , \notag \\
	= \  &\Big| \ \int_R \ \left(\frac{1}{\varphi'
	(\theta)}\frac{d}{d\theta}\right)^{N_2}\left( 
        \check\chi^{(-N_2)}_{(\frac{1}{4} , 4)}\big( \varphi(\theta) \big)
        \right)\cdot e^{-i\frac{(n-2)/2 + l}{\sqrt{r^2-m^2}} \ \theta} \
        d\theta \ \Big| \ , \notag \\
   \begin{split}\label{last_der_on_e_line}
	\lesssim \ &\int_R\ | \check\chi^{(-N_2)}_{(\frac{1}{4} , 4)}
	\big( \varphi(\theta) \big)|\cdot \left|\left(\frac{1}{\varphi'(\theta)}
	\frac{d}{d\theta}\right)^{N_2} 
	\left( e^{-i\frac{(n-2)/2 + l}{\sqrt{r^2-m^2}} \ \theta} 
	\right)\right|\ d\theta \\
	&\ + \ \ 
	\sum_{k=1}^{N_2}\ |
	\check\chi^{(-k)}_{(\frac{1}{4} , 4)}\big( \varphi(\theta) \big)|
	\cdot\big|\frac{1}{\varphi'(\theta)}
	\left(\frac{1}{\varphi'(\theta)}
	\frac{d}{d\theta}\right)^{k-1} 
	\left( e^{-i\frac{(n-2)/2 + l}{\sqrt{r^2-m^2}} \ \theta} 
	\right)\big|\ \Bigg|_{\partial R}
   \end{split}
\end{align}
To control the first term above, we bound:
\begin{align}
	&\left|\left(\frac{1}{\varphi'(\theta)}
	\frac{d}{d\theta}\right)^{N_2} 
	\left( e^{-i\frac{(n-2)/2 + l}{\sqrt{r^2-m^2}} \ \theta} 
	\right)\right|\ , \label{der_exp_bnd} \\
	\lesssim \ &\sup_{0\leqslant k\leqslant 2N_2-1}
	|\frac{1}{\varphi'(\theta)}|^k\cdot\sum_{j=0}^{N_2}\
	\left(
	\frac{1+l}{\sqrt{r^2-m^2}}\right)^j\, |P_{N_2-j}(\varphi^{(1+\bullet)}
	(\theta))|
	\ , \notag
\end{align}
where $P_0$ is a constant, and the other $P_k$ denote a homogeneous
expression of weight $k$ in the variables $(\varphi^{(2)},\ldots,
\varphi^{(i)},\ldots)$, where each $\varphi^{(i)}$ is given weight
$i-1$. Therefore, using the bound \eqref{phi_prime_bound} 
(which in particular implies that
$|1/\varphi'| \leqslant 2$) as well as the bounds
\eqref{phi_der_bnds}, we see that we can estimate:
\begin{align}
	\eqref{der_exp_bnd} \ \lesssim \ 
	&C_{N_2}\sum_{j=0}^{N_2}\ \left(\frac{1+l}{\sqrt{r^2-m^2}}\right)^j
	\cdot \frac{1}{(r-m)^{N_2-j}} \ , \notag \\
	\lesssim \ &C_{N_2}\left( \left(\frac{1+l}{\sqrt{r^2-m^2}}\right)^{N_2}
	+ \frac{1}{(r-m)^{N_2}}\right) \ . \notag
\end{align}
Finally, to deal with the boundary terms on the right hand side
of \eqref{last_der_on_e_line} we use the bounds \eqref{short_interval_h_bnd} 
and \eqref{chi_deriv_bnd} and simply estimate:
\begin{align}
	&\sum_{k=1}^{N_2}\ |
	\check\chi^{(-k)}_{(\frac{1}{4} , 4)}\big( \varphi(\theta) \big)|
	\cdot\Big|\frac{1}{\varphi'(\theta)}
	\left(\frac{1}{\varphi'(\theta)}
	\frac{d}{d\theta}\right)^{k-1} 
	\left( e^{-i\frac{(n-2)/2 + l}{\sqrt{r^2-m^2}} \ \theta} 
	\right)\Big|\ \Bigg|_{\partial R} \ , \notag\\
	\lesssim \ &C_{N_2}\sup_{0\leqslant k\leqslant N_2}
	|\check\chi^{(-k)}_{(\frac{1}{4} , 4)}\big( \varphi(\theta) \big)|
	\ \Bigg|_{\partial R} \ , \notag\\
	\lesssim \ &\frac{C_{N_1,N_2}}{(r-m)^{N_1}} \ . \notag
\end{align}\\

Combining the above estimates together, we have shown that:
\begin{equation}
        | I_1 | \ \lesssim \ \frac{C_{N_1,N_2}}{\sqrt{r^2-m^2}}
        \left( \frac{1}{(r - m)^{N_1}} \ + \ 
        \min_\pm \left\{  \left(\frac{l}{\sqrt{r^2-m^2}}\right)^{\pm N_2}  
        \right\}\right) \ . \notag
\end{equation}
This completes the proof of the asymptotic \eqref{psi_asym3}, and the
demonstration of Proposition \ref{psi_asym_prop}.
\end{proof}\ret

\ret

\section{Some $L^2$ dispersive estimates for the wave
equation; linear and bilinear estimates}

We are now ready to directly proceed with the proof of 
estimate \eqref{unit_ang_str_est}. As we have mentioned
previously, we may assume that $u_{1,N}$ is of the form
$e^{-it\sqrt{-\Delta}}\, f_{1,N}$, for some unit frequency and
dyadic angular frequency function $f_{1,N}$. For this  $u_{1,N}$, we 
the formulas from lines \eqref{f_wv_prop_exp}, \eqref{Hankel_wave_coef},
\eqref{c_fourier_series}, \eqref{c_wavelet_devel}, and \eqref{hankel_wavelet}
to expand its wave propagation into harmonics and Hankel-$\varphi$ transforms
as follows:
\begin{equation}
        e^{-it\sqrt{-\Delta}}\, f_{1,N}\, (r) \ = \ 
        \sum_{\substack{ l\sim N \\ i,k} } \ r^\frac{2-n}{2}\
        c^l_{i,k}\, \psi^l_{t - \frac{k}{4}}(r)\cdot Y^l_i \  
        \label{f1N_exp}
\end{equation}
In the above expression, we have absorbed the constants
$2\pi(\sqrt{-1})^l$ into each $c^l_{i,k}$. By orthogonality of everything
in sight, we have that:
\begin{equation}
        \lp{f_{1,N}}{L^2}^2 \ \sim \
        \sum_{\substack{i,k \\ l\sim N} }\ |c^l_{i,k}|^2 \ . \label{f1N_recovery_form}
\end{equation}
We now state our main result as follows:\\

\begin{prop}[$L^2$ dispersive estimate for the wave equation]
\label{L2_disp_est_prop}
For dimensions $2\leqslant n$, let $e^{-it\sqrt{-\Delta}}\, f_{1,N}$
be given by the formula \eqref{f1N_exp}. Then one has the following
estimate uniform in $t$ and $r$:
\begin{equation}
        \big| e^{-it\sqrt{-\Delta}}\, f_{1,N}\, (r) \big| \ \lesssim \
        \Big(\sum_{\substack{ l\sim N \\ i,k} } \ 
        \frac{|c^l_{i,k}|^2}{(1 + |t-\frac{4}{k}|)^{n-1}}\Big)^\frac{1}{2}
        \cdot N^\frac{n-1}{2} \ . \label{L2_disp_est}
\end{equation}
\end{prop}\ret

Before proceeding with the proof, let us first show briefly how the
above estimate may be used to show 
\eqref{unit_ang_str_est}--\eqref{n=2_unit_ang_str_est}. We begin with
\eqref{unit_ang_str_est}. What we 
need to do is to show that for every $0 < \eta$, there exists an
$\frac{2(n-1)}{n-2} < r_\eta$, such that the following estimate holds:
\begin{equation}
        \lp{e^{-it\sqrt{-\Delta}}\, f_{1,N} }{L^2(L^{r_\eta})} \ \lesssim \
        N^{\frac{1}{2} + \eta}\cdot\lp{f_{1,N}}{L^2} 
        \ , \label{sph_str_endpt_est}
\end{equation}
where the implicit constants depend on both $\eta$ and $r$. We also
need $r_\eta$ to approach $\frac{2(n-1)}{n-2}$ as $\eta\to 0$.
Now, interpolating\footnote{ This can be achieved by simply interpolating
in weighted $\ell^2$ and the usual Lebesgue
spaces using the map $\{c^l_{i,k}\} \to e^{-it\sqrt{-\Delta}}\, f_{1,N}$
for fixed time (see \cite{BLinterp}).}
\eqref{L2_disp_est} with the energy estimate
\eqref{basic_energy_est} gives the following result for $2\leqslant r_\eta$:
\begin{equation}
        \lp{e^{-it\sqrt{-\Delta}}\, f_{1,N} }{L^{r_\eta}_x} \ \lesssim \
        \Big(\sum_{\substack{ l\sim N \\ i,k} } \ 
        \frac{|c^l_{i,k}|^2}{ (1 + |t-\frac{4}{k}|)
        ^{ (\frac{1}{2}-\frac{1}{r_\eta})2(n-1)} }
        \Big)^\frac{1}{2}
        \cdot N^{ (\frac{1}{2} - \frac{1}{r_\eta})\,
        \frac{n-1}{2}} \ . \label{interp_disp_eng}
\end{equation}
Choosing $r_\eta$ according to the formula
$\frac{1}{2} + \eta = (\frac{1}{2} - \frac{1}{r_\eta})\,\frac{n-1}{2}$, we see that
we have the following identity holds:
$(\frac{1}{2}-\frac{1}{r_\eta})2(n-1) = 1 + 2\eta$. Therefore, for this
choice of $r_\eta$, we may square \eqref{interp_disp_eng} and integrate
directly in time to achieve \eqref{sph_str_endpt_est}. 
Also note that $r_\eta\to \frac{2(n-1)}{n-2}$ as $\eta\to 0$.\\

In the case of estimate \eqref{n=2_unit_ang_str_est}, that is $n=2$
spatial dimensions, we use H\"olders
inequality to see that
for every $0 < \eta$ the following estimate holds:
\begin{equation}
        \Big(\sum_{\substack{ l\sim N \\ i,k} } \ 
        \frac{|c^l_{i,k}|^2}{(1 +
        |t-\frac{4}{k}|)}\Big)^\frac{1}{2}
	\ \lesssim \ \Big(\sum_{\substack{ l\sim N \\ i,k} } \ 
        \frac{|c^l_{i,k}|^2}{(1 +
        |t-\frac{4}{k}|)^{1+\eta}}\Big)^\frac{1}{2(1+\eta)}\cdot
	\big(\sum_{\substack{ l\sim N \\ i,k} } \ |c^l_{i,k}|^2
        \big)^\frac{\eta}{2(1+\eta)} \ . \notag
\end{equation}
Therefore, of we choose $q_\eta = 2(1+\eta)$ a direct integration
and the energy bound \eqref{f1N_recovery_form} shows that we have:
\begin{equation}
        \lp{e^{-it\sqrt{-\Delta}}\, f_{1,N} }{L^{q_\eta}(L^\infty)} \ \lesssim \
        N^{\frac{1}{2}}\cdot\lp{f_{1,N}}{L^2} 
        \ . \notag
\end{equation}
We are now reduced to proving the estimate \eqref{L2_disp_est}.\\

\begin{proof}[proof of Proposition \ref{L2_disp_est_prop}]
Writing everything out, and using a Cauchy--Schwartz in the sum
over $i$ for each fixed $l$  we have that:
\begin{align}
        \big| e^{-it\sqrt{-\Delta}}\, f_{1,N}\, (r) \big| \ &\lesssim \
        \Big| \sum_{\substack{ l\sim N \\ i,k} } \ r^\frac{2-n}{2}\
        c^l_{i,k}\, \psi^l_{t - \frac{k}{4}}(r)\cdot Y^l_i
        \Big| \ , \label{L2_disp_start} \\
        &\lesssim \ \left( \sum_i
        \Big|\sum_{\substack{ l\sim N \\ k} } \ r^\frac{2-n}{2}\
        c^l_{i,k}\, \psi^l_{t - \frac{k}{4}}(r) \Big|^2 \right)^\frac{1}{2}
        \cdot \ \sup_{l\sim N}\, \left( \sum_i \ |Y^l_{i,k}|^2
        \right)^\frac{1}{2} \ . \notag
\end{align}
Now, using the formula in item (3) of Lemma \ref{basic_sph_lem}, and the fact
that $\frac{1}{l}(n+ 2l - 2)\binom{n+l-3}{l-1} \sim N^{n-2}$ for
$l\sim N$, we have:
\begin{align}
        \hbox{(L.H.S.)}\eqref{L2_disp_start} \ &\lesssim \
        \left( \sum_i
        \Big|\sum_{\substack{ l\sim N \\ k} } \ r^\frac{2-n}{2}\
        c^l_{i,k}\, \psi^l_{t - \frac{k}{4}}(r) \Big|^2 \right)^\frac{1}{2}
        \cdot \ N^\frac{n-2}{2} \ , \notag \\
        &\lesssim \ \left( \sum_{\substack{l\sim N \\ i}}
        \Big|\sum_{k } \ r^\frac{2-n}{2}\
        c^l_{i,k}\, \psi^l_{t - \frac{k}{4}}(r) \Big|^2 \right)^\frac{1}{2}
        \cdot \ N^\frac{n-1}{2} \ , \notag
\end{align}
where the last line above follows from a Cauchy--Schwartz and the bound
$l\sim N$. The proposition will now be shown if we can prove that:
\begin{equation}
        \sum_{k } \ r^\frac{2-n}{2}\
        |c^l_{i,k}|\cdot\big| \psi^l_{t - \frac{k}{4}}(r) \big| \ \lesssim \
        \Big(\sum_{\substack{ l\sim N \\ i,k} } \ 
        \frac{|c^l_{i,k}|^2}{(1 + |t-\frac{4}{k}|)^{n-1}}\Big)^\frac{1}{2}
        \ . \label{L2_disp_reduced}
\end{equation}
We now use the asymptotics of Proposition \ref{psi_asym_prop} and a
little computation to split:
\begin{equation}
        r^\frac{2-n}{2}\, \big| \psi^l_{t - \frac{k}{4}}(r) \big| \ \lesssim \
        A\Big(r,t - \frac{k}{4}\Big) \ + \ 
        B\Big(r,t - \frac{k}{4}\Big) \ , \notag
\end{equation}
where $B\Big(r,t - \frac{k}{4}\Big)$ is supported where 
$  1 < |t - \frac{k}{2}| < (r-1)   $, and:
\begin{align}
        A\Big(r,t - \frac{k}{4}\Big) \ &= \
        \frac{1}{(1 + |t - \frac{k}{4}|)^\frac{n-1}{2}}\cdot 
        \frac{1}{ (1 + \big| r - |t - \frac{k}{4}| \big|) } \ , \notag \\
        B\Big(r,t - \frac{k}{4}\Big) \ &= \
        \frac{1}{(1 + |t - \frac{k}{4}|)^\frac{n-1}{2}}\cdot
        \frac{1}{ (1 + \big| r - |t - \frac{k}{4}| \big| )^\frac{1}{2}}\cdot
        R\Big(l,t-\frac{k}{4},r\Big)
        \ . \notag 
\end{align}
Here $R\Big(l,t-\frac{k}{4},r\Big)$ is the function from line 
\eqref{psi_asym3_remainder}. Substituting $A\Big(r,t - \frac{k}{4}\Big)$
into the left hand side of the sum \eqref{L2_disp_reduced}, and using
a Cauchy--Schwartz immediately yields the desired result for this half
of things. In order to finish up, then, we only need substitute to
$B\Big(r,t - \frac{k}{4}\Big)$ into the right hand side. Doing
this and using a Cauchy--Schwartz and re-indexing, 
we see that it is enough to show:
\begin{equation}
        \sum_{\substack{ k \in \frac{1}{4}\cdot\mathbb{Z} + t \\
        |k| < (r-1)}} \ \frac{1}{(1 + \big| r - |k| \big|)}\cdot
        R(l,k,r) \ \lesssim \ 1 \ . \label{needed_k_sum}
\end{equation} 
uniform in $r$ and $l$. This can be done by comparing things to
the appropriate integrals and is left to the reader. Notice that
the convergence factor $R(l,k,r)$ avoids any logarithmic divergences
in this sum. This completes the proof of \eqref{L2_disp_est} and therefore 
the proof of Theorem \ref{full_ang_str_th}.
\end{proof}\ret

\ret

\section{Bilinear estimates for angularly regular data}

As is well known, for applications in the lower dimensional setting,
linear estimates of the form \eqref{full_ang_str_est} are not sufficient.
What is needed are multilinear versions. An extremely versatile
method for building these type of estimates is based on the ``fine and coarse
scale'' machine of T. Tao (see \cite{KRT}). The basic idea is to fix
a scale, say $\frac{1}{\mu}$ for $\mu\ll 1$, and then decompose the domain
of spatial variable into cubes with side lengths of this scale. Then,
one replaces the usual $L^r$ norm in the spatial variable with
$\ell^r(L^2)$, where the $L^2$ norm is taken on the ``fine'' scale
of each individual cube, while the $\ell^r$ norm represents the ``coarse''
scale which is summation over all cubes. One reason this method
is so important, is that it allows one to use the bilinear construction process
directly in an iteration procedure where  resorting to the canned estimates
that this method ultimately provides may be unduly burdensome. 
Also, the way these estimates are constructed will allow us to generate
bilinear estimates where only one term in the product contains extra
angular regularity. These type of estimates are extremely important
for applications when one proves inductive estimates via an
iteration procedure where it is necessary to consider estimates for
products with an angular derivative falling on one or the other
term in the product.
Because of these considerations,
we will content ourselves here with the dual scale estimates themselves, 
and not bother with listing out the various 
multilinear estimates which follow from them. In the $(4+1)$ and higher dimensions, 
the usual dual scale estimates read:\\

\begin{thm}[``Improved'' Strichartz estimates]\label{imp_str_th}
Let $4\leqslant n$ be the number of spatial dimensions. Let 
$0 < \mu \lesssim 1$ be given, and let $\{Q_\alpha\}$ be a partition 
of $\mathbb{R}^n$
into cubes of side length $\sim\frac{1}{\mu}$.
Then if 
$u_1$ is a unit frequency solution to the equation $\Box u_1=0$,
the following estimates hold:
\begin{equation}
        \lp{\big(\sum_\alpha \lp{u_1(t)}
        {L^2(Q_\alpha)}^r\big)^\frac{1}{r}}{L_t^2}
        \ \lesssim \ \mu^{-1} \, 
        \lp{u_1(0)}{L^2} 
        \ , \label{impr__str_est}
\end{equation}
where $\frac{2(n-1)}{n-3} \leqslant r$.
\end{thm}\ret 

\noindent We show here that both the range of \eqref{impr__str_est}
\emph{and} the power of $\mu$ that appearers there can be 
significantly improved. In what follows, we will only bother proving an
estimate which is the analog of \eqref{unit_ang_str_th}. More general
estimates which involve various amounts of angular regularity can
then be gained by interpolating this estimate with \eqref{impr__str_est}.
Also, it is not so easy to work out the Littlewood--Paley theory in the
angular variable for localized norms like those that appear on the
left hand side of \eqref{impr__str_est}. However, since we are already
loosing a small amount of angular regularity, one may simply replace the
Littlewood--Paley sum \eqref{ang_lp_sp_sum} by an $\ell^1$ sum. As we have noted
before, this has no bearing on applications.\\

\begin{thm}[``Improved'' frequency localized Strichartz estimates for 
angularly regular data; endpoint case]\label{imp_ang_str_th}
Let $3 \leqslant n$ be the number of spatial dimensions, 
and let $u_{1,N}$ be a unit frequency and angular frequency localized solution 
to the homogeneous  wave equation $\Box u_{1,N} = 0$. Let 
$0 < \mu \lesssim 1$ be given, and let $\{Q_\alpha\}$ be a partition 
of $\mathbb{R}^n$
into cubes of side length $\sim\frac{1}{\mu}$. Then for
every $0 < \eta$, there is a $C_\eta$ and $\frac{2(n-1)}{n-2} < r_\eta$ depending on $\eta$,
such that $r_\eta\to \frac{2(n-1)}{n-2}$ as $\eta\to 0$
such that the following estimate holds:
\begin{equation}
        \lp{\big(\sum_\alpha \lp{u_{1,N}(t)}
        {L^2(Q_\alpha)}^{r_\eta}\big)^\frac{1}{r_\eta}}{L_t^2}
        \ \lesssim \ C_\eta \, \mu^{-\frac{1}{2} -2\eta} \, 
        N^{\frac{1}{2}+\eta}\, \lp{u_{1,N}(0)}{L^2} 
        \ . \label{impr_ang_str_est}
\end{equation}
\end{thm}\ret

\begin{rem}
Up to the small loss in $\frac{1}{\mu}$ and $N$, the estimate \eqref{impr_ang_str_est}
is sharp when tested against the bilinear analog of the Knapp counterexamples
\eqref{wave_osc_int}. We construct these briefly as follows. We consider (frequency)
initial data sets $\chi_{B^\epsilon}$ which along with being highly localized in the
angular variable, are also well localized in the radial variable.
That is, we are now assuming $\chi_{B^\epsilon}$ is supported on a small
square of dimensions $\sim \epsilon\times\epsilon\times\ldots\times\epsilon$, lying along
the $\xi_1$ axis between $1/2 < \xi_1 < 2$.  A quick calculation then shows that the
integral:
\begin{equation}
        e^{it\sqrt{-\Delta}} f^\epsilon_1\, (x) \ = \
        \int e^{2\pi i (t |\xi| + x\cdot\xi)} \chi_{B^\epsilon}(\xi)\, d\xi
        \ , \notag
\end{equation}
behaves like $|e^{it\sqrt{-\Delta}} f^\epsilon_1\, (x)|
\sim \epsilon^{n}$ on the space--time region $S^\epsilon_{t,x}$:
\begin{align}
        t \ &= O(\epsilon^{-2}) \ , \notag \\
        t + x_1 \ &= \ O(\epsilon^{-1}) \ , \notag \\
        x' \ &= \ O(\epsilon^{-1}) \ . \notag
\end{align}
Choosing our cubes $Q_\alpha$ with side lengths $\frac{1}{\mu}
\sim\frac{1}{\epsilon}$, we see that for
any $\frac{2(n-1)}{n-2} < r$ (in fact for any $1\leqslant r \leqslant \infty$) the
following is true:
\begin{equation}
        \mu^{-\frac{1}{2}}\, \epsilon^{-\frac{1}{2}}\,
	\lp{f^\epsilon_1}{L^2}  \ \sim \ 
	\epsilon^n\cdot \frac{1}{\epsilon^\frac{n}{2}}\cdot
	\frac{1}{\epsilon} \ \lesssim \
        \lp{\big(\sum_\alpha \lp{  e^{it\sqrt{-\Delta}} f^\epsilon_1   }
        {L^2(Q_\alpha)}^{r}\big)^\frac{1}{r}}{L_t^2} 
\end{equation}
Therefore, using \eqref{Omega_of_f}, we see that \eqref{impr_ang_str_est} is indeed sharp
for this sequence of initial data. Of course the condition $\frac{2(n-1)}{n-2} < r$
cannot be improved, even for spherically symmetric initial data. In fact, one can 
also see for these type of waves (say with the asymptotic \eqref{good_sph_asym}), 
\eqref{impr_ang_str_est} is again sharp. We leave this simple calculation to the interested reader.  
\end{rem}\ret\ret

\begin{proof}[proof of estimate \eqref{impr_ang_str_est}]
The proof will be similar in spirit to that of
\eqref{unit_ang_str_est}. However, we will need to use orthogonality in a more
fundamental way here. This is because we will not be able to rely solely on 
$L^\infty$ as we did
in the proof of \eqref{L2_disp_est}. Ultimately, this has to do with the fact that
the wave packets $\psi^l_m(r)$ are not well localized in physical space when
$|m| < r$ and $r\ll l^2$ (see the remark after Proposition \ref{psi_asym_prop}).\\

In order to prove estimate \eqref{impr_ang_str_est}, we begin by writing the basic
energy estimate \eqref{basic_energy_est} in the following way (we are still using the
notation from line \eqref{f1N_exp}):
\begin{equation}
        \lp{\big(\sum_\alpha \lp{  e^{-it\sqrt{-\Delta}}f_{1,N}           }
        {L^2(Q_\alpha)}^2\big)^\frac{1}{2}}{L_t^\infty} \ \lesssim \
	\left( \sum_{\substack{l\sim N \\ i,k}} \ |c^l_{i,k}|^2\right)^\frac{1}{2}
        \ . \label{basic_energy_wp_form}
\end{equation}
Our next step is to prove the following spatially localized fixed time estimate:
\begin{equation}
        \sup_\alpha\ \lp{   e^{-it\sqrt{-\Delta}}f_{1,N}    }{L^2(Q_\alpha)} \ \lesssim \
	|\ln\mu|\mu^{-\frac{n-1}{2}} N^\frac{n-1}{2} \ \left(
	\sum_{\substack{l\sim N \\ i,k}} \ \frac{|c^l_{i,k}|^2}{(1 + |t-\frac{k}{4}|)^{n-1}}
	\right)^\frac{1}{2} \ . \label{imp_L2_disp_est}
\end{equation}
Assuming for the moment the validity of \eqref{imp_L2_disp_est}, we may interpolate
between \eqref{basic_energy_wp_form} and \eqref{imp_L2_disp_est} to achieve the 
following estimate for $2\leqslant r_\eta$:
\begin{multline}
        \left( \sum_\alpha \
        \lp{e^{-it\sqrt{-\Delta}}\, f_{1,N}}{L^2(Q_\alpha)}^{r_\eta}\right)
	^\frac{1}{r_\eta}\\
	\lesssim \ \left( \sum_{\substack{ l\sim N \\ i,k} } \
        \frac{|c^l_{i,k}|^2}{ (1 + |t-\frac{4}{k}|)^{ (\frac{1}{2}-\frac{1}{r_\eta})(n-1)} }
        \right)^\frac{1}{2}
        \cdot |\ln\mu| \, \mu^{-(\frac{1}{2} -\frac{1}{r_\eta})\,
        \frac{n-1}{2}}   N^{ (\frac{1}{2} -\frac{1}{r_\eta})\,
        \frac{n-1}{2}} \ \label{imp_interp_est}
\end{multline}
As in the previous section, choosing  $r_\eta$ by the identity,
$\frac{1}{2} + \eta = (\frac{1}{2} - \frac{1}{r_\eta})\,\frac{n-1}{2}$, we have that
$(\frac{1}{2}-\frac{1}{r_\eta})2(n-1) = 1 + 2\eta$. Therefore, squaring \eqref{imp_interp_est}
and integrating directly in time, we will have achieved \eqref{impr_ang_str_est}.
Therefore, we now concentrate on proving \eqref{imp_L2_disp_est}.\\

It suffices to show \eqref{imp_L2_disp_est} for a fixed $\alpha$. Therefore, we will now
assume that we are on a fixed cube $Q_\alpha$. Our first step is foliate  $Q_\alpha$
with the hypersurfaces $\mathbb{S}_r^{n-1}\cap Q_\alpha$, where $\mathbb{S}_r^{n-1}$
is the sphere of radius $r$ centered at the origin. Notice that for each $r$, one
has the area estimate $|\,\mathbb{S}_r^{n-1}\cap Q_\alpha| \lesssim \mu^{-(n-1)}$.
Therefore, using a Cauchy--Schwartz, it suffices to prove the estimate:
\begin{equation}
        \lp{   e^{-it\sqrt{-\Delta}}f_{1,N}    }{L_r^2(L^\infty(\mathbb{S}_r^{n-1}))
	  (Q_\alpha)} \ \lesssim \
	|\ln\mu| \, N^\frac{n-1}{2} \ \left(
	\sum_{\substack{l\sim N \\ i,k}} \ \frac{|c^l_{i,k}|^2}{(1 + |t-\frac{k}{4}|)^{n-1}}
	\right)^\frac{1}{2} \ . \label{imp_L2_disp_est2}
\end{equation}
Next, we chop $Q_\alpha$ into at most $|\ln\mu|$ dyadic pieces, which are of the form
$R_j\cap Q_\alpha$, where $R_j=\{ x \ \big| \ 2^j < |x| < 2^{j+1}\}$ is the radial
dyadic region of size $r\sim 2^j$. We only need to do this for $0 < j$, i.e we keep the
ball of bounded radius (say radius $=2$) as a single region $R_0$. Therefore, to show
\eqref{imp_L2_disp_est2}, it suffices to estimate:
\begin{equation}
        \sup_j\ \lp{   e^{-it\sqrt{-\Delta}}f_{1,N}    }{L_r^2(L^\infty(\mathbb{S}_r^{n-1}))
	  (R_j)} \ \lesssim \ N^\frac{n-1}{2} \ \left(
	\sum_{\substack{l\sim N \\ i,k}} \ \frac{|c^l_{i,k}|^2}{(1 + |t-\frac{k}{4}|)^{n-1}}
	\right)^\frac{1}{2} \ . \label{imp_L2_disp_est3}
\end{equation}
Now fix $R_j$. To prove estimate \eqref{imp_L2_disp_est3} here,
we run the decomposition \eqref{f1N_exp} and decompose the resulting sum into the 
sum of pieces:
\begin{equation}
        e^{-it\sqrt{-\Delta}}f_{1,N} \ = \ \Sigma_0 + \Sigma_1 \ , \notag 
\end{equation}
where:
\begin{align}
        \Sigma_0 \ &= \ \sum_{\substack{ l\sim N \\ i,k \\ |t-\frac{4}{k}| > 4\cdot 2^j
	} } \ r^\frac{2-n}{2}\
        c^l_{i,k}\, \psi^l_{t - \frac{k}{4}}(r)\cdot Y^l_i \ , \notag \\
	\Sigma_1 \ &= \ \sum_{\substack{ l\sim N \\ i,k \\ |t-\frac{4}{k}| \leqslant 4\cdot 2^j
	} } \ r^\frac{2-n}{2}\
        c^l_{i,k}\, \psi^l_{t - \frac{k}{4}}(r)\cdot Y^l_i \ , \notag
\end{align}
To estimate \eqref{imp_L2_disp_est3} on the sum $\Sigma_1$, we simply keep it as a whole
object and use the Bernstein inequality \eqref{Rn_sph_bernstein} to estimate:
\begin{align}
        &\sup_j\ \lp{   \Sigma_1   }{L_r^2(L^\infty(\mathbb{S}_r^{n-1}))
	  (R_j)} \ , \notag \\
	\lesssim \ &N^\frac{n-1}{2}\, 2^{-\frac{n-1}{2}j}\, \lp{ \Sigma_1 }{L_x^2}
	\ , \notag \\
	\lesssim \ &N^\frac{n-1}{2}\, 2^{-\frac{n-1}{2}j}\,
	\left( \sum_{\substack{ l\sim N \\ i,k \\ |t-\frac{4}{k}| \leqslant 4\cdot 2^j
	} } \
        |c^l_{i,k}|^2\right) \ , \notag\\
	\lesssim \ &N^\frac{n-1}{2}\left( 
	\sum_{\substack{ l\sim N \\ i,k \\ |t-\frac{4}{k}| \leqslant 4\cdot 2^j} } \
        \frac{|c^l_{i,k}|^2}{(1 + |t-\frac{k}{4}|)^{n-1}} \right) \ . \notag
\end{align}
As was to be shown. Therefore, we are reduced to bounding $\Sigma_0$. Notice that
for each fixed $r\in R_j$, and for each term in this sum, we have 
$|t-\frac{k}{4}| > r$. Therefore, we can use the well localized asymptotic 
\eqref{psi_asym1}
and a Cauchy--Schwartz in the $k$ summation (in a way similar to the 
computation started on line
\eqref{L2_disp_start}) to bound this term in $L^\infty(\mathbb{S}_r^{n-1}))$ by:
\begin{equation}
        |\Sigma_0 (r)| \ \lesssim \
	\left( \sum_{\substack{ l\sim N \\ i,k \\ |t-\frac{4}{k}| > 4\cdot 2^j
	} } \  \frac{|c^l_{i,k}|^2}{(1 + |t-\frac{k}{4}|)^{n-1}}
	\cdot\frac{1}{(1 + \big| r - |t-\frac{t}{4}| \big|)^2} \right)^\frac{1}{2}
	\cdot N^\frac{n-1}{2}\ . \notag 
\end{equation}
Squaring this last expression, and integrating each term in the sum with respect to 
$dr$ (notice there is no extra volume element because 
we took the sup on $\mathbb{S}_r^{n-1}$), we see that we have proved 
\eqref{imp_L2_disp_est3} for this portion of things. This then completes the
proof of estimate \eqref{imp_L2_disp_est}, and therefore the proof of 
\eqref{impr_ang_str_est}.
\end{proof}

\ret

\section{Appendix}

We present here a simplified proof of Theorem \ref{unit_ang_str_th} and 
Theorem \ref{imp_ang_str_th}.
 First, we note that in order
to prove both \eqref{unit_ang_str_est} 
\begin{equation} \label{eq:ang-str}
\|u_{1,N} \|_{L^2_t L^{p_\eta}_x} \les C_\eta N^{\frac 12+\eta} \|u_{1,N}(0)\|_{L^2}
\end{equation}
and \eqref{impr_ang_str_est}
 \begin{equation} \label{eq:ang-bilin}
\|u_{1,N} \|_{L^2_t \ell^{p_\eta} L^{2}_\mu} \les C_\eta \mu^{-\frac 12-\eta}
N^{\frac 12+\eta} \|u_{1,N}(0)\|_{L^2}
\end{equation}
with $p_\eta\to 2(n-1)/(n-2-\eta)$ and arbitrary small $\eta>0$, 
it suffices to prove the estimate \eqref{impr_ang_str_est}. Here $u_{1,N}$ 
is a unit frequency solution of the wave equation $\Box u=0$ of angular frequency 
$N$ and the norm
$$
\|f\|_{\ell^{p} L^{2}_\mu} = \Big (\sum_{\alpha} \|f\|_{L^2(Q_\alpha)}^p\Big )^{\frac 1p},
$$
where $\{Q_\alpha\}$ is a partition of ${\Bbb R}^n$ into cubes $Q_\alpha$ of
side length of $\mu^{-1}$.

To see this, we choose a partition $\{Q_\alpha\}$ of size $\mu=1$
and compute, using the Sobolev embedding on $Q_\alpha$, that for any $2\le p$:
\begin{align}
	\lp{u_{1,N} }{L^2_tL^p_x} \ &= \ 
	\lp{\left(\sum_{\alpha} \lp{u_{1,N}
	}{L^p(Q_\alpha)}^p\right)^\frac{1}{p}}{L^2_t} \ , \notag\\
	&\lesssim\ \lp{\left(\sum_{\alpha} \lp{u_{1,N}\,
	}{L^2(Q_\alpha)}^p\right)^\frac{1}{p}}{L^2_t} \ . \notag
\end{align}
Therefore, we see that it suffices to deal with the estimate
\eqref{eq:ang-bilin}. First, using the 
Sobolev inequality $L^\infty_\omega\subset \langle\Omega \rangle ^{-\frac {n-1}2} L^2_\omega$ 
on the unit sphere ${\Bbb S}^{n-1}$ for angular frequency localized
functions,
 we have the following estimate
for any tiling of $\mathbb{R}^n$ by cubes $\{Q_\alpha\}$ of side length
$\sim\mu^{-1}$:
\begin{align}
	\lp{u_{1,N}(t)}{L^2(Q_\alpha)}^2 \ &\les 
	\int_{r\omega\in Q_\alpha} |{u_{1,N}(t)}|^2 dx\,\les\, 
	\mu^{-\frac {n-1}2} \int_0^\infty dr\, r^{n-1} \sup_{\omega\in {\Bbb S}^{n-1}}
	|u_{1,N} (r\omega,t)|^2 \nonumber\\
	&\lesssim  \
	\mu^{-\frac{n-1}{2}} \ \lp{
	(1+r)^{-\frac{n-1}{2}}\, \langle\Omega\rangle^{\frac{n-1}{2}} u_{1,N}(t)  }{L^2_x}
	\ . \nonumber
\end{align} 
Thus,
\begin{equation}
\|u_{1,N}(t)\|_{\ell^\infty L^2_\mu}\, \les\,\, 
\mu^{-\frac{n-1}{2}} N^{\frac {n-1}2} \ \lp{
	(1+r)^{-\frac{n-1}{2}}\,  u_{1,N}(t)  }{L^2_x}.\nonumber
\end{equation}	
Interpolating this with the trivial  estimate
$$
\|u_{1,N}(t)\|_{\ell^2 L^2_\mu} \,\les\, 
 \ \lp{
	\,  u_{1,N}(t)  }{L^2_x},
$$
 we arrive at the following
estimate which will be our point of departure:
\begin{equation}
	\|u_{1,N}\|_{L^2_t \ell^{p_\eta} L^2_\mu}
	\ \lesssim \ C_\eta \ \mu^{-\frac{1}{2} - \eta} 
	\lp{ (1+r)^{-\frac{1}{2}(1 + \eta)}\langle \Omega^{\frac{1}{2}(1 + \eta)}\rangle u_{1,N}}{L_t^2L_x^2} \ , \label{igors_reduction} 
\end{equation}
where $p_\eta =\frac {2(n-1)}{n-2-\eta}$. Therefore, using \eqref{igors_reduction}
and that $\langle\Omega\rangle^{\frac{1}{2}(1 + \eta)} u_{1,N}$ is another unit
frequency solution of the wave equation localized at the angular frequency $N$,
we see that in order to prove \eqref{eq:ang-bilin}, it suffices
to prove the following space-time Morawetz type estimate for unit frequency solutions
of the homogeneous wave equation when $3\leqslant n$:
\begin{equation}
	\int_{-\infty}^\infty \int_{\mathbb{R}^n}\
	\frac{1}{(1+r)^{1 + \eta}} |u_1(t,x)|^2 \ dx dt \ \lesssim \
	\lp{u_1(0)}{L_x^2}^2 \ . \label{comm_est}
\end{equation}
In the above estimate, the implicit constant depends on $\eta$ and the 
dimension, and from now on we will keep this dependence implicit.
Estimate \eqref{comm_est}, as well as as its version with $\eta=0$ and a 
logarithmic loss in time, was proved in the work of Keel-Smith-Sogge
\cite{KSS} in dimension three with the help of the sharp Huygen's 
principle.
We will prove \eqref{comm_est} directly using essentially nothing but
an integration by part argument.  Before we continue, let us make one
more reduction. It turns out that \eqref{comm_est} is more naturally
proved for the spatial gradient $\nabla_x u_1$. Since $u_1$ is unit frequency,
this reduction does not effect the validity of \eqref{comm_est}.
Indeed, suppose we have the following estimate:
\begin{equation}
	\int_{-\infty}^\infty \int_{\mathbb{R}^n}\
	\frac{1}{(1+r)^{1 + \eta}} |\nabla_x u_1(t,x)|^2 \ dx dt \ \lesssim \
	\lp{\nabla_{t,x} u_1\, (0)}{L_x^2}^2 \ . \label{grad_comm_est}
\end{equation}
for any unit frequency solution $u_1$. Let $\tilde u_1$ be another unit
frequency solution with the property that $\Delta \tilde u_1=-u_1$. 
Then for any $\epsilon>0$,
\begin{align*}
\int_{-\infty}^\infty &\int_{\mathbb{R}^n}\
	\frac{1}{(1+r)^{1 + \eta}} |u_1(t,x)|^2 \ dx dt \ = 
\int_{-\infty}^\infty \int_{\mathbb{R}^n}\
	\frac{1}{(1+r)^{1 + \eta}} \Delta \tilde u_1 \overline {u}_1 \ dx dt\\
	&=\int_{-\infty}^\infty \int_{\mathbb{R}^n}\
	\frac{1}{(1+r)^{1 + \eta}} \nabla_x \tilde u_1 \nabla_x\overline {u}_1 \ dx dt	+(1+\eta) \int_{-\infty}^\infty \int_{\mathbb{R}^n}\
	\frac{r \nabla_x r}{(1+r)^{3 + \eta}} \nabla_x \tilde u_1 \overline{u}_1 \ dx dt\\
&	\lesssim \epsilon^{-1}\int_{-\infty}^\infty \int_{\mathbb{R}^n}\
	\frac{1}{(1+r)^{1 + \eta}} |\nabla_x \tilde u_1(t,x)|^2  dx dt +
	\int_{-\infty}^\infty \int_{\mathbb{R}^n}\
	\frac{1}{(1+r)^{1 + \eta}} |{\nabla_x u_1(t,x)}|^2 \ dx dt \\ 
	&+
	\epsilon\int_{-\infty}^\infty \int_{\mathbb{R}^n}\
	\frac{1}{(1+r)^{1 + \eta}} |u_1(t,x)|^2 \ dx dt 
\end{align*}
and \eqref{comm_est} follows from \eqref{grad_comm_est} by choosing 
a sufficiently small $\epsilon$.
Therefore, we may now assume that we are trying to achieve estimate
\eqref{grad_comm_est}. The
fact that \eqref{grad_comm_est} contains the gradient $\nabla_x$ will
allow us to prove it for an arbitrary solution to the
homogeneous wave equation. This will be done 
using some more or less standard energy-momentum tensor techniques. 
Let $\phi$ be a solution to the homogeneous wave equation:
\begin{equation}
        \Box \phi \ = \ 0 \ ,
\end{equation}
and let $Q_{\alpha\beta}[\phi]$ be its energy-momentum tensor:
\begin{equation}
        Q_{\alpha\beta}[\phi] \ = \ \partial_\alpha \phi\,
        \partial_\beta\phi  - \frac{1}{2}g_{\alpha\beta}
	\partial^\gamma\phi\partial_\gamma\phi \ . \notag
\end{equation}
Here the Greek indices run on the set $\alpha,\beta,\gamma = 0, \ldots
, n$, and $g_{\alpha\beta} = \hbox{diag}(-1,1,\ldots,1)$ is the
standard Minkowski metric. The key feature of $Q_{\alpha\beta}[\phi]$
is that it is space-time divergence free:
\begin{equation}
        D^\alpha Q_{\alpha\beta}[\phi] \ = \ 0 \ ,
\end{equation}
where $\nabla$ denotes the Levi-Civita connection of 
$g_{\alpha\beta}$. We now contract $Q_{\alpha\beta}[\phi]$
a radial vector-field:
\begin{equation}
        X \ = \ f(r)\partial_r \ . \notag
\end{equation} 
Define the momentum density:
\begin{equation}
        P_\alpha [\phi,X] \ = \  Q_{\alpha\beta}[\phi] X^\beta \ . \notag
\end{equation}
A quick calculation shows that the divergence of $P_\alpha [\phi,{}X]$
satisfies the identity:
\begin{equation}
        D^\alpha P_\alpha [\phi,{}X] \ = \ \frac{1}{2}
	Q_{\alpha\beta}[\phi] \, \pi^{\alpha\beta} \ , \label{div_calc}
\end{equation}
where:
\begin{equation}
        \pi_{\alpha\beta} \ = \ D_\alpha X_\beta
	+ D_\beta X_\alpha\ , \notag
\end{equation}
is the deformation tensor of $X$. Introducing the orthonormal
frame $\{\partial_t , \partial_r, e_A\}$, where the $\{e_A\}$ form
an orthonormal frame on each $n-1$ sphere $S_{t,r} =
\{r=\hbox{const.}\}$, one can easily calculate this quantity to be:
\begin{equation}
       \pi_{\alpha\beta} \ = 
   \begin{pmatrix} 
	0 & 0 & \ldots & 0 \\
	0 & f'(r)  & \ldots & 0 \\
	0 & 0 & \ \ \frac{f(r)}{r} \delta_{AB}  
   \end{pmatrix} \ , \label{def_calc}
\end{equation}
where $\delta_{AB}$ denotes the metric on the spheres $S_{t,r}$.
In particular, we see that:
\begin{equation}
        \hbox{tr}\pi  \ = \
	f'(r) + (n-1) \frac{f(r)}{r} \ , \label{tr_pi_form}
\end{equation}
Next, using \eqref{div_calc} and \eqref{def_calc}, a direct computation shows that:
\begin{equation}
        2 \, D^\alpha P_\alpha [\phi, X] \ = \
	f'(r) |\partial_r\phi|^2 + \frac{f(r)}{r} |\snabla \phi|^2
	-\frac{1}{2} \hbox{tr}\pi 
	\partial^\gamma\phi\,\partial_\gamma\phi \ . \notag
\end{equation}
In the above formula, $|\snabla \phi|^2 = \delta^{AB}|e_A(\phi)|\cdot
|e_B(\phi)|$ denotes the angular portion of the spatial 
gradient $|\nabla_x\phi|^2$.  Since, 
$
\partial^\gamma\phi\,\partial_\gamma\phi =\frac 12 \Box |\phi|^2$,
if  we define the modified momentum density:
\begin{equation}
        \widetilde{P}_\alpha [\phi,{}X] \ = \ 
	P_\alpha [\phi,{}X] + \frac{1}{4} \hbox{tr}{}\pi 
	\phi\,\partial_\alpha\phi - \frac{1}{8}
	\partial_{\alpha}\left(\hbox{tr}{}\pi\right)
	|\phi|^2 \ , \notag
\end{equation}
we end up with the identity:
\begin{equation}
        D^\alpha \widetilde{P}_\alpha [\phi,{}X] \ = \
	\frac{1}{2}\left(f'(r) |\partial_r\phi|^2 + \frac{f(r)}{r} |\snabla
	\phi|^2\right) -\frac{1}{8}
	\Delta\left( \hbox{tr} {}\pi 
	\right)\, |\phi|^2\ . \label{mod_div}
\end{equation}
Integrating \eqref{mod_div} over a time slab, we arrive at the
following a-priori estimate for $\phi$:
\begin{multline}
        \int_{\mathbb{R}^n} \ \widetilde{P}_0 [\phi,{}X](0)\ dx \ \
	\ 
	= \ \ \ \int_{\mathbb{R}^n} \ \widetilde{P}_0 [\phi,{}X
	](T)\ dx
	\ \ \ + \ \ \ \\
        \int_0^T\int_{\mathbb{R}^n}\ \frac{1}{2}\left(f'(r)
	|\partial_r\phi|^2 
	+ \frac{f(r)}{r} |\snabla\phi|^2\right) -\frac{1}{8}
	\Delta\left( \hbox{tr} {}\pi
	\right)\, |\phi|^2 \ dx\, dt \ , \label{basic_apriori}
\end{multline}
where:
\begin{equation}
        \int_{\mathbb{R}^n} \ \widetilde{P}_0 [\phi,{}X](0)\ dx
	\ = \ \int_{\mathbb{R}^n}\Big (X^\alpha\, \partial_t\phi(0) \, \partial_\alpha\phi(0)
	+ \frac{1}{4} \hbox{tr} {} \pi 
	\phi(0)\partial_t\phi(0) \Big )\ dx \ . \label{boundary_terms}
\end{equation}
with an identical expression for the time $=T$ boundary piece on the
right hand side of \eqref{basic_apriori} above. In particular, using 
\eqref{tr_pi_form}, one has that if $|f(r)| \lesssim 1$ then:
\begin{align}
        \left|\int_{\mathbb{R}^n} \ \widetilde{P}_0 
	[\phi,{}X](0)\ dx \right| \ &\lesssim \
	\lp{ r^{-1} \phi\, (0)}{L^2_x}\cdot
	\lp{  \nabla_{t,x}\phi\, (0)}{L^2_x} \ , \notag\\
	&\lesssim \ \lp{\nabla_{t,x}\phi\, (0)}{L^2_x}^2 \ , \label{boundary_energy_est}
\end{align}
where in the last line we used the Hardy inequality.
Similar estimate holds for the other boundary term 
$\widetilde{P}_0 [\phi,{}X](T)$. Therefore,
using the usual conservation of energy, we have the following 
a-priori estimate for solutions to the homogeneous wave equation:
\begin{multline}
        \big| \ \int_{-\infty}^{\infty}\int_{\mathbb{R}^n}\ \frac{1}{2}\left(f'(r)
	|\partial_r\phi|^2 
	+ \frac{f(r)}{r} |\snabla\phi|^2\right) -\frac{1}{8}
	\Delta\left( \hbox{tr}{}\pi 
	\right)\, |\phi|^2 \ dx\, dt\ \big| \\ \lesssim \ \ \
	\ \lp{\nabla_{t,x}\phi\, (0)}{L^2_x}^2 \ . \label{basic_apriori1}
\end{multline}
We now use \eqref{basic_apriori1} to derive \eqref{grad_comm_est} by choosing
the weight function $f(r)$. 
$$
f(r) = \frac {r}{\epsilon +r}
$$
for some $\epsilon >0$. A direct calculation shows that, with this choice 
of $f$, we have 
$$
\Delta \left( \hbox{tr}{}\pi 
	\right) = 
	\Delta \big (\frac {\epsilon}{(\epsilon +r)^2} + \frac {n-1}{\epsilon+r}\big )=
	-\frac 1{r(\epsilon +r)^3} \Big ((n-3) r +3(n-3)\frac {\epsilon r}{\epsilon+r}
	+\frac {3\epsilon^2(n-1)}{\epsilon+r} \Big ) < 0
$$
in dimensions $n\ge 3$.
Therefore this gives the a-priori estimate:
\begin{equation}
        \int_{-\infty}^{\infty}\int_{\mathbb{R}^n}\ 
	\Big (\frac{\epsilon}{(\epsilon+r)^2} |\partial_r\phi|^2+ 
	\frac{1}{\epsilon+r} |\snabla\phi|^2 \Big )
	\ dx\, dt \ \lesssim \ \lp{\nabla_{t,x}\phi\, (0)}{L^2_x}^2
	\  \label{basic_mor1}
\end{equation}
with an implicit constant in $\lesssim$ independent of $\epsilon$. 
In particular, choosing $\epsilon = 2^k$ with an integer $k\ge 0$, we 
obtain 
$$
  \int_{-\infty}^{\infty}\int_{|x|\le 1}\ 
	{|\nabla_x\phi|^2}
	\ dx\, dt   +     \int_{-\infty}^{\infty}\int_{2^{k-1}\le |x|\le 2^{k+1}}\ 
	\frac {|\nabla_x\phi|^2}{r}
	\ dx\, dt \ \lesssim \ \lp{\nabla_{t,x}\phi\, (0)}{L^2_x}^2.
$$
Dividing the above inequality by $2^{-\eta k}$ and summing over $k$
immediately yields the desired estimate \eqref{grad_comm_est}.
\ret\ret


\end{document}